\documentclass[12pt]{amsart}

\title[Diamond for $\mathfrak{sl}(n)$]{Diamond representations of $\mathfrak{sl}(n)$}
\author[D. Arnal]{Didier Arnal}
\address{
Institut de Math\'ematiques de Bourgogne\\
UMR CNRS 5584\\
Universit\'e de Bourgogne\\
U.F.R. Sciences et Techniques
B.P. 47870\\
F-21078 Dijon Cedex\\France} \email{Didier.Arnal@u-bourgogne.fr}
\author[N. Bel Baraka]{Nadia Bel Baraka}
\address{
Institut de Math\'ematiques de Bourgogne\\
UMR CNRS 5584\\
Universit\'e de Bourgogne\\
U.F.R. Sciences et Techniques
B.P. 47870\\
F-21078 Dijon Cedex\\France} \email{Nadia.Bel-Baraka@u-bourgogne.fr}
\author[N. Wildberger]{Norman J. Wildberger}
\address{
School of Mathematics\\
University of New South Wales\\
Sydney 2052\\
Australia} \email{n.wildberger@unsw.edu.au}

\date{07/11/05}

\usepackage{graphicx}
\usepackage{rotate,epic,eepic,epsf,epsfig}

\usepackage{amsmath}
\usepackage{amsfonts}
\usepackage{amssymb}
\newtheorem{defi}{Definition}
\newtheorem{prop}{Proposition}
\newtheorem{theo}{Theorem}
\newtheorem{lem}{Lemma}
\newtheorem{rema}{Remark}

\begin{document}
\maketitle \vspace{.20cm}
\noindent
{\bf Abstract}\\

In \cite{W}, there is a graphic description of any irreducible,
finite dimensional $\mathfrak{sl}(3)$ module. This construction,
called diamond representation is very simple and can be easily
extended to the space of irreducible finite
dimensional ${\mathcal U}_q(\mathfrak{sl}(3))$-modules.\\

In the present work, we generalize this construction to
$\mathfrak{sl}(n)$. We show this is in fact a description of the
reduced shape algebra, a quotient of the shape algebra of
$\mathfrak{sl}(n)$. The basis used in \cite{W} is thus naturally
parametrized with the so called quasi standard Young tableaux. To
compute the matrix coefficients of the representation in this
basis, it is possible to use Groebner basis for the ideal of
reduced Pl\"ucker relations defining the reduced shape algebra.

\medskip
\section{Introduction}

\

In this paper, we consider the irreducible finite dimensional
representations of the Lie algebra $\mathfrak{sl}(n) =
\mathfrak{sl}(n,{\mathbb C})$. Of course these representations are
well known and there are very explicit descriptions for them, for
instance in \cite{FH}.\\

First, $\mathfrak{sl}(n)$ acts naturally on ${\mathbb C}^n$, its
fundamental representations are the natural actions on ${\mathbb
C}^n,\wedge^2{\mathbb C}^n,\dots,\wedge^{n-1}{\mathbb C}^n$, they
have highest weights $\omega_1,\dots,$ $\omega_{n-1}$. Each simple
$\mathfrak{sl}(n) $-module has a highest weight $\lambda$ and this
highest weight characterizes the module. Note ${\mathbb
S}^\lambda$ this module,
 it is a submodule of the tensor product
$$Sym^{a_1}({\mathbb C}^n)\otimes Sym^{a_2}(\wedge^2{\mathbb C}^n)\otimes
\dots\otimes Sym^{a_{n-1}}(\wedge^{n-1}{\mathbb C}^n),$$
if $\lambda=a_1\omega_1+\dots+a_{n-1}\omega_{n-1}$.\\

The direct sum ${\mathbb S}^{\bullet}$ of all the simple modules
has a natural realization as the shape algebra of
$\mathfrak{sl}(n)$, {\sl i.e.} as the algebra ${\mathbb
C}[SL(n)]^{N^+}$ of polynomial functions on the group $SL(n)$,
which are invariant under the right multiplication by upper
triangular matrices. Let $g$ be an element in $SL(n)$, denote
$\delta^{(s)}_{i_1,\dots,i_s}(g)$ the determinant of the submatrix
of $g$ obtained by considering the $s$ first columns of $g$ and
the rows $i_1<\dots<i_s$, then ${\mathbb S}^{\bullet}$ is
generated as an algebra by the functions
$\delta^{(s)}_{i_1,\dots,i_s}$. More precisely, it is the quotient
of ${\mathbb C}[\delta^{(s)}_{i_1,\dots,i_s}]$ by the ideal
$P(\delta)$ generated by the Pl\"ucker relations.\\

Generally a parametrization of a basis for ${\mathbb S}^\lambda$ is given by the
set of semi-standard Young tableaux $T$ of shape $\lambda$ {\sl i.e.} with
$a_{n-1}$ columns of size $n-1$, $\dots$, $a_1$ columns of size 1.\\

Using this description, we give here a natural ordering on the set
of variables $\delta^{(s)}_{i_1,\dots,i_s}$, we determine the
Groebner basis of $P(\delta)$ for this ordering, getting the
corresponding basis of the  quotient as monomials
$\delta^T$, for $T$ semi-standard.\\

Thus the action of upper triangular matrices on this basis can be easily
computed. (See for instance the description given in \cite{LT}).\\

On the other hand, in \cite{W}, N. Wildberger gave a really
different presentation of the simple $\mathfrak{sl}(3)$-modules.
This description is based on the construction of the diamond cone
for $\mathfrak{sl}(3)$, it is an infinite dimensional
indecomposable module for the Heisenberg Lie algebra with a very
explicit basis. The matrix coefficients are integral numbers and
fixing the highest weight $\lambda$, it is easy to build the
corresponding representation of $\mathfrak{sl}(3)$, on the
submodule generated by this vector in the diamond
cone.\\

In this paper, we extend this presentation to $\mathfrak{sl}(n)$.
In fact the diamond cone module  is a quotient of the shape
algebra. We call this quotient the reduced shape algebra. It is
the quotient of ${\mathbb C}[\delta^{(s)}_{i_1,\dots,i_s}]$ by the
ideal $P_{red}(\delta)$ sum of the ideal of Pl\"ucker relations
and the
ideal generated by $\delta^{(s)}_{1,\dots,s}-1$.\\

With the same approach as above, we define a new ordering on the variables
$\delta^{(s)}_{i_1,\dots,i_s}$, with this ordering, we can compute the Groebner
basis for $P_{red}(\delta)$ and the corresponding basis for the quotient : the
set of monomials $\delta^T$, for  some Young tableaux $T$ called here
quasi-standard. The action of the upper triangular matrices on this basis is
easy to compute : this gives us the diamond cone for $\mathfrak{sl}(n)$.\\

In order to refind the complete $\mathfrak{sl}(n)$-modules, we
have to define a symmetry on each $\mathbb{S}^{\lambda}$ and on
the corresponding submodule in the reduced shape algebra. This
symmetry exchanges the role of $N^+$ and $N^-$ and we get the
complete $\frak{sl}(n)$ representation.\\

Unfortunately, this symmetry corresponds to a modification of
the ordering on Young tableaux, thus, if $n>3$ to a different
basis in $\mathbb{S}^{\lambda}$. The $\frak{n}^-$ action on the first
base is not so simple as in \cite{W}.\\

\medskip
\section{Usual (algebraic) presentation of the $\mathfrak{sl}(n)$ simple
modules}

\

Let us consider the Lie algebra $\mathfrak{sl}(n) =
\mathfrak{sl}(n, \mathbb C)$: it is the set of $n\times n$
traceless matrices, it is the Lie algebra of the Lie group $SL(n)$
of $n\times n$ matrices, with determinant 1. The Cartan algebra
$\mathfrak{h}$ is the space of diagonal matrices:
$$\mathfrak{h}=\left\{H=\left[\begin{matrix}\theta_1&&0\cr&\ddots&\cr0&&\theta_n
\end{matrix}\right],\quad \theta_j\in{\mathbb C},\quad\theta_1+\dots+\theta_n=0
\right\}.$$

We put $\alpha_i(H)=\theta_i$. The root system of $\mathfrak{sl}(n)$ is the set
of linear form on $\mathfrak{h}$ generated by the $\alpha_i-\alpha_j$, ($i\neq
j$).

The usual basis $\Delta$ for the root system is given by :
$$\Delta=\{ \alpha_i-\alpha_{i+1},~~i=1, 2,\dots,n-1\}$$
The root space corresponding to the simple root $\eta_i=\alpha_i-\alpha_{i+1}$
is generated by the upper triangular matrix:
$$X_{\eta}=\left[\begin{array}{clccc}0&&&&\\&&1&&\\&&\ddots&&\\
&&&&0\end{array}\right].$$
The root space corresponding to $-\eta$
is generated by lower triangular matrix:
$$Y_{\eta} = \left[\begin{array}{crccc}0&&&&\\&&\ddots&&\\&1&&&\\
&&&&0\end{array}\right] = ~^tX_{\eta}$$ these matrices generate
$\mathfrak{sl}(n)$ as a Lie algebra.

\

A weight $\lambda$ for $\mathfrak{sl}(n)$ is a linear form :
$$\lambda:\left[\begin{matrix}\theta_1&&0\cr&\ddots&\cr0&&\theta_n
\end{matrix}\right]\mapsto(a_1+\dots+a_{n-1})\theta_1+(a_2+\dots+a_{n-1})
\theta_2+\dots+a_{n-1}\theta_{n-1}.$$ If $a_1,\dots,a_{n-1}$ are
positive integral numbers, we shall say that $\lambda$ is a
dominant integral weight. This is the case if and only if
$\lambda$ is a linear combination $\lambda = \sum^{n-1}\limits_{j
= 1} a_j \omega_j$, with positive integral coefficients $a_j$, of
the fundamental weights:
$$\omega_j=\alpha_1+\dots+\alpha_j:\left[\begin{matrix}\theta_1&&0\cr&\ddots&\cr
0&&\theta_n\end{matrix}\right]\mapsto\theta_1+\dots+\theta_j\qquad(1\leq j\leq
n-1).$$

\

The set of simple $\mathfrak{sl}(n)$-modules up to equivalence is isomorphic to
the set of dominant integral weights. More precisely, $\mathfrak{sl}(n)$ acts
naturally on $V={\mathbb C}^n$ (with canonical basis $e_1,\dots,e_n$), thus also
on the totally antisymmetric tensor products $\wedge^jV$ ($j=1,\dots,n-1$) and
on the symmetric tensor products $Sym^{a_j}(\wedge^jV)$ and finally on
$$Sym^{a_1}(V)\otimes Sym^{a_2}(\wedge^2V)\otimes\dots\otimes Sym^{a_{n-1}}
(\wedge^{n-1}V).$$

For each dominant integral weight $\lambda=\sum a_j\omega_j$, the corresponding
simple module ${\mathbb S}^\lambda(V)$ is the submodule of
$$Sym^{a_1}(V)\otimes Sym^{a_2}(\wedge^2V)\otimes\dots\otimes Sym^{a_{n-1}}
(\wedge^{n-1}V).$$
generated by the vector:
$$v^\lambda=(e_1)^{a_1}\otimes(e_1\wedge e_2)^{a_2}\otimes\dots\otimes
(e_1\wedge\dots\wedge e_{n-1})^{a_{n-1}}.$$

With this construction, we get each simple $\mathfrak{sl}(n)$-module, and two
distinct weights $\lambda$, $\lambda'$ give rise to inequivalent simple
$\mathfrak{sl}(n)$-modules.\\

This action gives rise by exponentiation to a representation of $SL(n)$. Let us
put
$$\Omega=\left[\begin{matrix}0&&&&\varepsilon_n\cr&&&.&\cr
&&.&&\cr&.&&\cr\varepsilon_n&&&&0\end{matrix}\right]$$
where $\varepsilon_n=1$ if $\left[\frac{n}{2}\right]$ is even and $\varepsilon_n
=e^{\frac{i\pi}{n}}$ if $\left[\frac{n}{2}\right]$ is odd. Then $\Omega$
belongs to $SL(n)$. In fact, this matrix, acting by adjoint action generates the
longest element of the Weyl group of $SL(n)$. It correspnds to a change in the
choice of simple roots and nilpotent subalgebras $\frak{n}^+$ and $\frak{n}^-$,
if $X=\left[x_{ij}\right]$ is a strictly upper triangular matrix, $\Omega^{-1}X
\Omega=\left[x_{(n+1-i)(n+1-j)}\right]$ is strictly lower triangular. Let us
put:
$$v^\lambda_-=(e_n)^{a_1}\otimes(e_n\wedge e_{n-1})^{a_2}\otimes\dots\otimes(e_n
\wedge\dots\wedge e_2)^{a_{n-1}}=\varepsilon_n^{-|\lambda|}\Omega.v^\lambda,$$
with $|\lambda|=a_1+2a_2+\dots+(n-1)a_{n-1}$. Then $v^\lambda_-$ is a lowest
weight vector in $\mathbb{S}^\lambda(V)$.

\

\section{The shape algebra: abstract algebraic presentation}

\

Let us put:
$$
{\mathbb S}^\bullet(V)=\bigoplus_{\lambda}~~{\mathbb S}^\lambda(V).
$$
Since we have an explicit realization of each highest weight vector, it is
possible to define a natural comultiplication $\Delta$ on ${\mathbb S}^\bullet
(V)$, just by defining
$$
\Delta:{\mathbb S}^\lambda(V)\longrightarrow\bigoplus_{\mu+\nu=\lambda}
~~{\mathbb S}^\mu(V)\otimes{\mathbb S}^\nu(V)
$$
as the unique $\mathfrak{sl}(n)$-morphism sending $v^\lambda$ on
$$
\Delta(v^\lambda)=\sum_{\mu+\nu=\lambda}v^\mu\otimes v^\nu.
$$
$\Delta$ is cocommutative. The contragredient module $({\mathbb
S}^{\lambda})^*$ is naturally identified with ${\mathbb
S}^{^t\lambda}$ where $^t\lambda = \sum a_{n-i}\omega_i~~~{\rm
if}~~~ \lambda = \sum a_i \omega_i$. By transposition, $\Delta$
defines a commutative multiplication $m$ on ${\mathbb
S}^\bullet(V)$:
$$m = {^t\Delta}:{\mathbb S}^{^t\mu}(V)\otimes{\mathbb S}^{^t\nu}(V)\longrightarrow
{\mathbb S}^{^t\mu+^t\nu}(V).$$
By definition, if $\mu=\sum_jb_j\omega_j$, $\nu=\sum_jc_j\omega_j$,
$$m(v^\mu\otimes v^\nu)=v^\mu.v^\nu=v^{\mu+\nu}=e_1^{b_1+c_1}\otimes
\dots\otimes(e_1\wedge\dots\wedge e_{n-1})^{b_{n-1}+c_{n-1}}.$$

Since each isotypic component of the $SL(n)$ module ${\mathbb
S}^\bullet(V)$ is simple the multiplication $m$ is characterized
by this relation and the condition
$$m\left({\mathbb S}^\mu(V)\otimes{\mathbb S}^\nu(V)\right)\subset
{\mathbb S}^{\mu+\nu}(V).$$

\

We shall call \underbar{shape algebra} of $SL(n)$ the algebra
${\mathbb S}^\bullet(V)$ equipped with the above multiplication.

\

\

By construction the shape algebra is generated as an algebra by
the subspace $V\oplus\wedge^2V\oplus\dots\oplus\wedge^{n-1}V$.
Thus it is a quotient of the algebra denoted in \cite{FH}:
\begin{align*}
A^\bullet(V)&=Sym^\bullet\left(V\oplus\wedge^2V\oplus\dots\oplus\wedge^{n-1}V
\right)\\
&=\bigoplus_{a_1,\dots,a_{n-1}}Sym^{a_{n-1}}(\wedge^{n-1}V)\otimes\dots\otimes
Sym^{a_1}(V).
\end{align*}

\

We define now the ideal of Pl\"ucker relations: it is the ideal $P$ of
$A^\bullet(V)$ generated by the vectors in $Sym^2(\wedge^pV)$:
\begin{align*}
(e_{i_1}\wedge\dots\wedge e_{i_p}).&(e_{j_1}\wedge\dots\wedge e_{j_p})+\\
&+\sum_{\ell=1}^p(-1)^\ell(e_{j_1}\wedge e_{i_1}\wedge\dots\wedge\widehat{
e_{i_\ell}}\wedge\dots\wedge e_{i_p}).(e_{i_\ell}\wedge e_{j_2}\wedge\dots\wedge
e_{j_p})
\end{align*}
and by the vectors in $\wedge^pV\otimes\wedge^qV$ ($p>q$)
\begin{align*}
(e_{i_1}\wedge\dots\wedge e_{i_p}).&(e_{j_1}\wedge\dots\wedge e_{j_q})+\\
&+\sum_{\ell=1}^p(-1)^\ell(e_{j_1}\wedge e_{i_1}\wedge\dots\wedge\widehat{
e_{i_\ell}}\wedge\dots\wedge e_{i_p}).(e_{i_\ell}\wedge e_{j_2}\wedge\dots\wedge
e_{j_q}).
\end{align*}

\

\begin{theo}

\

The shape algebra ${\mathbb S}^\bullet(V)$ is the quotient of $A^\bullet(V)$ by
the ideal $P$.
\end{theo}

This theorem is well known. There is a complete proof in
\cite{FH} p. 235, this result is cited by Towber in \cite{LT}
as a theorem due to Kostant.

\

We define a symmetry $\tau$ in $\mathbb{S}^{\bullet}$ just by
putting:
$$\tau(v)=\Omega.v  ~~~~~~{\rm if}~~~~v \in \mathbb{S}^{\bullet}(V)$$
Since the multiplication is a morphism of $\frak{sl}(n)$ and
$SL(n)$ module, $\tau(vv')=\tau(v)\tau(v')$. Especially, we can
define the multiplication just as above by fixing $v^\lambda_-.v^\mu_-=
v^{(\lambda + \mu)}_-$.

\

Now for each matrix $A$ in $\frak{sl}(n)$, $\Omega A \Omega=
{^\tau A}$ is the matrix defined by a central symmetry on the
entries of $A$:
$$ {^\tau A} = [a_{n+1-i,n+1-j}] ~~~{\rm if}~~~ A=[a_{i,j}]$$
If $A$ is the matrix $X_\eta$ for a positive root $\eta=\alpha_i -
\alpha_j$, ${^\tau X_\eta} = \Omega X_\eta\Omega$ is the matrix
$Y_{\tau_\eta}$ if $^\tau\eta$ is the positive root $^\tau\eta =
\alpha_{n+1-j}-\alpha_{n+1-i}$ Then for each $v$ in
$\mathbb{S}^\bullet$:
$$(\tau\circ X_\eta \circ\tau)(v) = \Omega X_\eta \Omega~ v = Y_{\tau_\eta} v$$

\

\section{The shape algebra: geometric presentation}

\

The shape algebra can also be viewed as an algebra of functions on
a quotient ${SL(n)}/{N^+}$ of the Lie group $SL(n)$. Denote $N^+$
the subgroup of matrices $n^+ = \left[\begin{matrix}1&&&*\\
&& \ddots &&\\
0&&&1
\end{matrix}\right]$ \\

Let us consider the space ${\mathbb C}[SL(n)]={\mathbb C}[g_{ij}]/
(det-1)$ of all polynomial functions $f$ with respect to the
entries $g_{ij}$ of the matrix $g\in SL(n)$. There is a
$SL(n)\times SL(n)$ action on this space, defined as follows:
$$((g_1,g_2).f)(g')=f(^tg_1g'g_2).$$

Since this space is generated by the invariant finite dimensional subspaces of
class of functions with degree less than $N$ ($N=0,1\dots$), this action is
completely reducible in a sum of finite dimensional simple $SL(n)\times SL(n)$
modules. The highest vector for these modules are class of functions $f$ such
that:
$$f(^tn_1^+gn_2^+)=f(g)~~~~~~~~n^+_1 \in N^+,~n^+_2 \in N^+.$$
But, let us consider the restriction of $f$ to the dense set of
$g$ such that, for $s = 1, \dots, n$,
$\delta^{(s)}_{1,2,\dots,s}(g)\neq0$. On this set, using the Gauss
method, we can reduce $g$ to a diagonal matrix, getting:
$$g= {^tn}_1^+\left[\begin{matrix}\delta_1^{(1)}(g)&&&&0\\&
\displaystyle\frac{\delta_{1,2}^{(2)}(g)}{\delta_1^{(1)}(g)}&&&\\&&\ddots&&\\
&&&\displaystyle\frac{\delta_{1,2,\dots,n-1}^{(n-1)}(g)}
{\delta_{1,\dots,n-2}^{(n-2)}(g)}&\\0&&&&\displaystyle\frac{1}
{\delta_{1,2,\dots,n-1}^{(n-1)}(g)}\end{matrix}\right]n_2^+.$$
The highest weight is $(\lambda,\lambda)$ with $\lambda=\sum a_i\omega_i$, that
means $f$ is a polynomial function in the variables
$$\delta_1^{(1)}(g),\frac{\delta_{1,2}^{(2)}(g)}{\delta_1^{(1)}(g)},\dots,
\frac{\delta_{1,2,\dots,n-1}^{(n-1)}(g)}{\delta_{1,\dots,n-2}^{(n-2)}(g)},
\frac{1}{\delta_{1,2,\dots,n-1}^{(n-1)}(g)},$$
homogeneous with degree $a_1+\dots+a_{n-1},a_2+\dots+a_{n-1},\dots,
a_{n-1},0$, {\sl i.e.} the function $f$ is a multiple of the function:
$$\delta^\lambda=\left(\delta_1^{(1)}\right)^{a_1}\left(\delta_{1,2}^{(2)}
\right)^{a_2}\dots\left(\delta_{1,2,\dots,n-1}^{(n-1)}\right)^{a_{n-1}}.$$
Acting with only the first factor $SL(n)$ on these functions, we get all the
polynomial $N^+$ right invariant functions on $SL(n)$. Due to the form of the
bi-invariant functions $f$, these functions are polynomial functions in the
$\delta$-variables :
$${\mathbb C}[SL(n)]^{N^+}\simeq{\mathbb C}[\delta^{(s)}_{i_1,\dots,i_s}]/
P(\delta),$$
where $P(\delta)$ is an ideal. Moreover each irreducible representation of
$SL(n)$ happens exactly one times in this space, thus as a vector space,
$${\mathbb S}^\bullet({\mathbb C}^n)\simeq{\mathbb C}[SL(n)]^{N^+}.$$\\

Acting on $\delta^\lambda$ ($\lambda=\sum_ia_i\omega_i$) on the
left by $N^- = {^t(N^+)}$, we get polynomial functions which
contains only monomials of the form:
$$\prod_{s=1}^{n-1}\prod_{k=1}^{a_s}\delta^{(s)}_{i_1^k,\dots,i_s^k}.$$
Let us call $V^{a_1,\dots,a_{n-1}}$ the space of such functions. In view of our
description, it is a simple module and the isotypic component of type $\lambda$
in ${\mathbb C}[SL(n)]^{N^+}$.

Finally the usual pointwise multiplication of polynomial functions send
$V^{a_1,\dots,a_{n-1}}\otimes V^{b_1,\dots,b_{n-1}}$ into $V^{(a_1+b_1),\dots,
(a_{n-1}+b_{n-1})}$. Thus the above identification
$${\mathbb S}^\bullet({\mathbb C}^n)\simeq{\mathbb C}[SL(n)]^{N^+},$$
characterized by $v^\lambda\mapsto\delta^\lambda$ is a morphism of algebra.
\

\begin{prop}

\

The shape algebra is isomorphic to the algebra ${\mathcal O}(SL(n)/N^+)$ of the
regular functions on the homogeneous space $SL(n)/N^+$.

The ideal $P(\delta)$ is the ideal generated by the Pl\"ucker relations written
on the $\delta$ functions.

\end{prop}

\

\begin{rema}

\

In this presentation of ${\mathbb S}^{\bullet}({\mathbb C}^n)$, the $SL(n)$
action on the elements of the shape algebra, viewed as a polynomial function $f$
is very natural since it is just:
$$(g.f)(g') = f(^tgg'), ~~ g \in SL(n), ~g' \in SL(n). $$
\end{rema}

\

The symmetry $\tau$ can be directly implemented in the space ${\mathbb
C}[SL(n)]^{N^+}$. Indeed $\tau$ is up to conjugation by $\Omega$, a morphism of
$SL(n)$ modules and the formula
$$\tau(e_1\wedge\dots\wedge e_s)=\varepsilon_n^se_n\wedge\dots\wedge e_{n+1-s}$$
becomes here
$$\tau(\delta^{(s)}_{1,2,\dots,s})=\varepsilon_n^s\delta^{(s)}_{n,(n-1),\dots,
(n+1-s)}.$$
But, if we put for any regular function $f$ on $SL(n)$, $(\theta f)(g)=f(\Omega
g)$, we define a bijection from ${\mathbb C}[SL(n)]^{N^+}$ into itself such that
$$g.\theta(f)=\theta(\Omega^{-1}g\Omega.f)\quad \hbox{and}\quad\theta(
\delta^{(s)}_{1,\dots,s})=\delta^{(s)}_{n,n-1,\dots,n+1-s}.$$
Thus $\tau=\theta$.\\

\

\section{The shape algebra : Combinatorial presentation}

\

The usual basis of $ S^\lambda(V)$ are parameterized by the semi standard Young
tableaux with shape $\lambda$. Let us be more precize:

\

We can naturally associate to each $\delta$ variable a column $C$:
$$\delta^C=\delta^{(p)}_{i_1,\dots,i_p}\longrightarrow\begin{array}{|c|}
\hline
i_1 \\
\hline
i_2 \\
\hline
\vdots \\
\hline
i_p \\
\hline
\end{array}$$
Then if we identify two Young tableaux which differ only by a permutation of
their columns, the set of Young tableaux defines a linear basis for the algebra
$\mathbb{C}[\delta^{(s)}_{i_1,\dots,i_s}]$:

$$\delta^T=\delta^{(p_1)}_{i_1,\dots,i_{p_1}}\delta^{(p_2)}_{j_1,\dots,j_{p_2}}
\dots\delta^{(p_k)}_{\ell_1,\dots,\ell_{p_k}}\longrightarrow
\begin{array}{l}
\begin{array}{|c|c c|c c|}
  \hline
   & \dots &  &  &  \\
  \hline
\end{array}\\
\begin{array}{|c|ccc|}
   & \vdots &  &    \\
   \hline
\end{array}\\
\begin{array}{|c|c|}
   & \dots    \\
     & \\
   \hline
   \end{array}\\
\begin{array}{|c|}
 \\
\hline
\end{array}
\end{array}$$
$(p_1\leq p_2\leq\dots\leq p_k)$. That means, we read the Young tableau
from right to left, using the following convention: if two different columns $C$
and $C'$ have the same height, we put in the first place in $T$ the column $C$
if
$$i_p=i'_p,i_{p-1}=i'_{p-1},\dots,i_{r+1}=i'_{r+1},~\hbox{and}~i_r<i'_r.$$

The Pl\"ucker relations are quadratic in the $\delta$ variables, they correspond
to linear combination of Young tableaux with two columns, for instance, we get
for $\mathfrak{sl}(3)$ the following relation between tableaux:
$$\delta^{(2)}_{12} \delta^{(1)}_{3} -  \delta^{(2)}_{23} \delta^{(1)}_{1}
+  \delta^{(2)}_{13} \delta^{(1)}_{2} \longrightarrow
\begin{array}{l}
\begin{array}{|c|c|}
\hline 1 & 3\\
\hline
\end{array}\\
\begin{array}{|c|}
\hline 2\\
\hline \end{array} \end{array} - \begin{array}{l}
\begin{array}{|c|c|}
\hline 2 & 1\\
\hline
\end{array}\\
\begin{array}{|c|}
\hline 3\\
\hline \end{array} \end{array} + \begin{array}{l}
\begin{array}{|c|c|}
\hline 1 & 2\\
\hline
\end{array}\\
\begin{array}{|c|}
\hline 3\\
\hline \end{array} \end{array}$$
In order to describe a basis for the quotient space:
$$ \mathbb{S}^{\bullet}(V) = \mathbb{C}[\delta^{(j)}_{i_1,\dots, i_j}]/_{\hbox{
$P(\delta)$}},$$
we will use the notion of Groebner basis \cite{CLO}.\\

Let us consider the algebra $\mathbb{C}[X_1, \dots , X_k]$ of polynomials in
the variable $X_i$ and an ideal $I$ of $\mathbb{C}[X_1, \dots , X_k]$.

Suppose we fix an ordering on the set of monomials $X^{\alpha_1}_1\dots
X^{\alpha_k}_k$ in $\mathbb{C}[X_1,$ $\dots,X_k]$ (for instance by using the
lexicographie ordering on words $\alpha_1\dots\alpha_k$ if we put
$X_1< X_2<\dots<X_k$). Then any polynomial $g$ has an unique leading term $LT(g)
$; the greatest monomial happening in $g$ for this ordering.\\

\begin{defi}

\

A finite subset $\{g_1,\dots,g_k\}$ of an ideal $I$ is said to be a reduced
Groebner basis for $I$ if and only if the leading term of any element of $I$ is
divisible by one of the leading term of $g_i$ and if for all $g_i$ no monomial
of $g_i$ is divisible by the leading term of some $g_j$ $j\neq i$.
\end{defi}

\

If $\{g_1,\dots,g_k\}$ is a reduced Groebner basis for $I$, then the set of
(classes of) monomials which are not divisible by any monomials $LT(g_i)$ ($i=1,
\dots,k$) is a basis of the quotient $\mathbb{C}[X_1, \dots , X_k]/I$.\\

Following \cite{FH}, we know there is in the ideal $P(\delta)$ the following
elements for any $p \geq q \geq r$:
$$
\delta^{(p)}_{i_1,i_2,\dots,i_p}\delta^{(q)}_{j_1,j_2,\dots,j_q} +
\sum\limits_{\begin{array}{c}
               A \subset\{i_1,\dots,i_p\}\\
               \# A = r
               \end{array}}\pm
\delta^{(p)}_{(\{i_1,\dots,i_p\}\setminus A)\cup\{j_1,\dots,j_r\}}
\delta^{(p)}_{A\cup \{j_{r+1},\dots,j_q\}}\eqno{(*)}$$
where $\delta^{(p)}_{(\{i_1,\dots,i_p\}\setminus A)\cup
\{j_1,\dots,j_r\}} = 0$ if there is a repetition of some index
and, if $\{k_1,\dots,k_p\}=(\{i_1,\dots,i_p\}\setminus A)\cup \{j_1,\dots,j_r\}$
and $k_1 <\dots<k_p$, then
$$\delta^{(p)}_{(\{i_1,\dots,i_p\}\setminus A)\cup\{j_1,\dots,j_r\}}=
\delta^{(p)}_{k_1,\dots, k_p}.$$

Now we put an ordering on the variables $\delta^{(p)}_{i_1,\dots,i_p}$ by the
following relations:
$$ \delta^{(1)}_{\dots}>\delta^{(2)}_{\dots}>\delta^{(n-1)}_{\dots}$$
and $\delta^{(p)}_{i_1,\dots,i_p}>\delta^{(p)}_{j_1,\dots,j_p}$ if $i_p=j_p,
\dots, i_{r+1}=j_{r+1}$ and $i_r<j_r$.

We put the lexicographic ordering on the monomials $\delta^T$ in $\mathbb{C}[
\delta^{(p)}_{i_1,\dots,i_p}]$.\\

\begin{rema}

\

In \cite{FH}, an ordering $<$ on Young tableaux is defined,in fact our
ordering is the reverse ordering since:
$$\delta^T < \delta^{T'}~~\hbox{if and only if}~~T' < T.$$
\end{rema}

\

Recall that a Young tableau is semi standard if its entries are increasing along
each row (and strictly increasing along each column). It is well known that the
set of semi standard Young tableau gives a basis for $\mathbb{C}[
\delta^{(p)}_{i_1,\dots,i_p}]/P(\delta)$ (see \cite{FH} for instance).\\

Our  ordering defines an unique Groebner basis for $P(\delta)$. We shall now
build this basis.

For any non semi standard Young tableau $T$ with 2 columns, there exists an
element in $P(\delta)$ of the form $(*)$. This relation can be written as:
$$\delta^T+\sum\limits_{j=1}^{n}\pm\delta^{T_j}~~\hbox{with}~~\delta^{T_j}<
\delta^{T}~~\forall j.$$
Each $T_j$ has the same shape as $T$ but some of them can be non semi standard.
We repeat the construction for each non semi standard $T_j$ and finally we get,
for each non semi standard $T$ with 2 columns, an element $f_T$ in $P(\delta)$
such that the leading term of $f_T$ is $\delta^T$ and all the monomials of $f_T$
have the form $a.\delta^{T'}$ with $T'$ semi standard and $\delta^{T'}<\delta^T$.

\begin{theo}

\

The set
$$G = \left\{f_S,~~S ~{\rm non ~semi~standard~with~2
            ~columns}\right\}$$
is the reduced Groebner basis of $P(\delta)$ for our ordering.
\end{theo}

\

\noindent
{\bf Proof:}

First denote $NS$ the set  of all monomials $\delta^T$ with $T$ non semi
standard. Since each non semi standard $T$ has 2 consecutive columns such that
the sub tableau defined by these 2 columns is non semi standard, $\delta^T$ is
divisible by one of the $\delta^S$, i.e. by one of the leading term of $f_S$.

Thus the ideal $<\delta^S>$ generated by the leading terms of $G$ contains the
vector space $span(NS)$.

Conversely let $T$ be a semi standard Young tableau. Suppose $T$ belongs to the
ideal $<LT(P(\delta))>$ generated by the leading terms of all the $f$ in
$P(\delta)$. That means:
$$\delta^T=f-\sum\limits_{T'<T} a_{T'}\delta^{T'}.$$
If any $T'$ is semi standard we keep this relation. If some of the $T'$ are non
semi standard, then $\delta^{T'}$ is in $<\delta^S>$ thus in $<LT(P(\delta))>$
and we repeat the construction for $\delta^{T'}$. We get finally:
$$\delta^T=f_0-\sum_{\begin{smallmatrix}
            T''< T\\
            T''{\rm semi~standard}
            \end{smallmatrix}} a_{T''}\delta^{T''}, ~~f_0 \in P(\delta).$$
This implies that
$$\delta^T+\sum\limits_{T''}a_{T''}\delta^{T''}\in P(\delta).$$
But this is impossible, since the set $\{\delta^T, T~~{\rm semi~standard}\}$ is a
basis for $\mathbb{C}[\delta^{(p)}_{i_1,\dots,i_p}]/P(\delta)$ Thus:
$$<LT(P(\delta))> = span(NS)$$
Moreover, since any monomial in $f_S$ is either $\delta^S$ or $a_T\delta^T$ with
$T$ semi standard, it can not be divisible by a $\delta^{S'}$ with $S'\neq S$,
$S'$ non semi standard with two columns. This proves our theorem.\\

The usual basis of the shape algebra $\mathbb{S}^\bullet(V)$ by semi standard
Young tableaux can thus be described as a natural basis of a the quotient of the
polynomial algebra $\mathbb{C}[\delta^{(p)}_{i_1,\dots,i_p}]$ by the ideal of
Pl\"ucker relations, if we put the ordering $<$ on the monomials $\delta^T$.\\

Especially, we can write the action of any element of the Lie algebra
$\mathfrak{sl}(n)$ on any polynomial function with variables
$\delta^{(s)}_{i_1,\dots,i_p}$, for instance, if $X_\alpha=E_{ij}$ $i\neq j$
then $X_\alpha$ acts on $\mathbb{C}[\delta^{(p)}_{i_1,\dots,i_p}]$ as the
derivation:
$$X_\alpha f = \frac{d}{ds}|_{s=0}f(exp~s^t X_\alpha.)=
\sum\limits_{\{i_1,\dots,i_p\}\cap\{i,j\}=\{j\}}\pm\delta^{(p)}_{(\{i_1,\dots,
i_p\}\setminus\{j\})\cup\{i\}}\frac{\partial f}{\partial\delta^{(p)}_{i_1,\dots,
i_p}}.$$
Finally, the Cartan algebra acts on $f$ as the derivation
$$H f = \frac{d}{ds}|_{s=0}f(exp~s^t H.)=
\sum(\theta_{i_1}+\dots+\theta_{i_p})\delta^{(p)}_{i_1,\dots,i_p}\frac{\partial
f}{\partial\delta^{(p)}_{i_1,\dots,i_p}},$$
if $H=\left[\begin{matrix}\theta_1&&0\cr&\ddots&\cr0&&\theta_n
\end{matrix}\right]$. This action defines the action on the quotient by
$P(\delta)$, since we have a Groebner basis for the ideal $P(\delta)$, the
quotient action on the basis of semi standard Young tableaux reduces to compute
the canonical form of the polynomial $X_\alpha f$ or $Hf$, this is easy to
do with usual computer software.

As an illustration, we gives a graphic description of the $N^+$ part of the
adjoint representation ${\mathbb S}^{\omega_1+\omega_2}({\mathbb C}^3)$ of
$\mathfrak{sl}(3)$ (see \cite{K} for similar presentation).

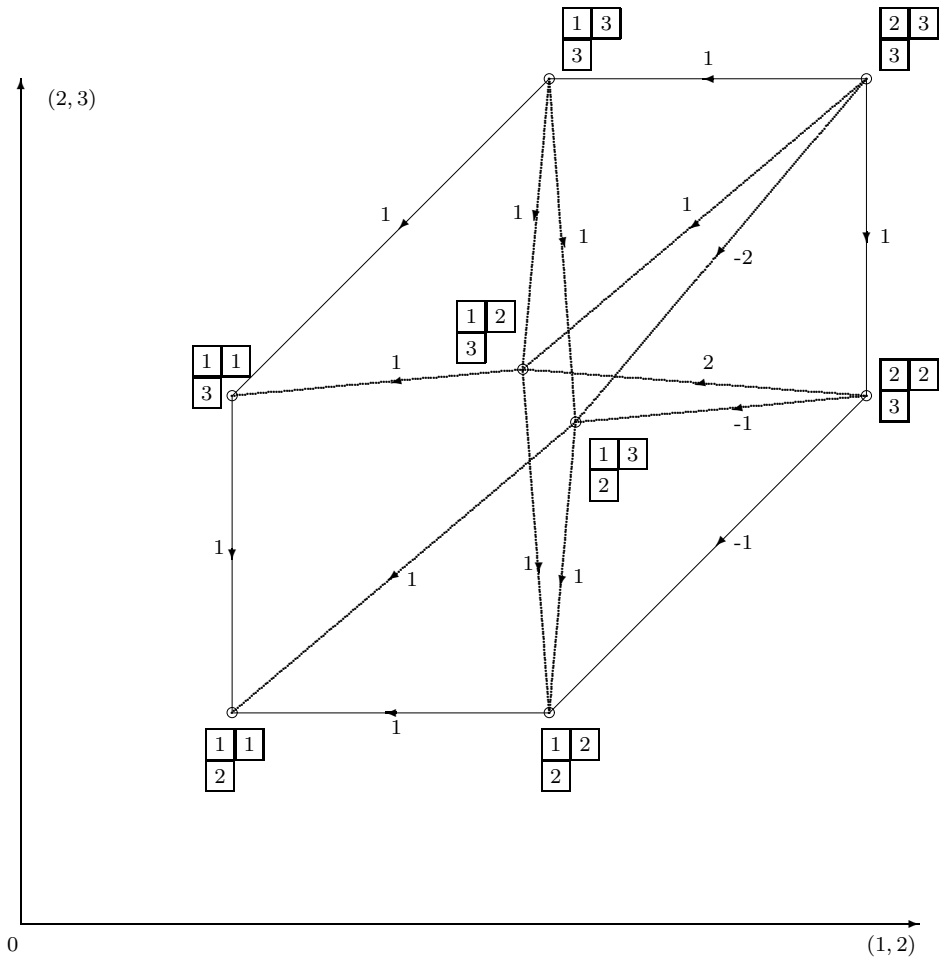
\begin{center}

{\tiny
\begin{picture}(410,600)(0,0)
\put(10,40){\vector(1,0){340}} 
\put(10,40){\vector(0,1){320}}
\put(20,350){$(2,3)$} 
\put(330,30){$(1,2)$}
\put(5,30){0}
\path(90,120)(210,120)
\path(90,120)(90,240)
\path(210,120)(330,240)
\path(330,240)(330,360)
\path(330,360)(210,360)
\path(210,360)(90,240)
\dottedline[.](330,360)(200,250)
\dottedline[.](330,360)(220,230)
\dottedline[.](210,360)(200,250)
\dottedline[.](210,360)(220,230)
\dottedline[.](330,240)(200,250)
\dottedline[.](330,240)(220,230)
\dottedline[.](90,240)(200,250)
\dottedline[.](210,120)(200,250)
\dottedline[.](90,120)(220,230)
\dottedline[.](210,120)(220,230)
\put(90,120){\circle{4}}
\put(90,240){\circle{4}}
\put(210,120){\circle{4}}
\put(200,250){\circle{4}}
\put(220,230){\circle{4}}
\put(330,240){\circle{4}}
\put(330,360){\circle{4}}
\put(210,360){\circle{4}} 
\put(202,100){$\begin{array}{l}
\framebox{1}\framebox{2}\\
\framebox{2}
\end{array}$}
\put(75,100){$\begin{array}{l}
\framebox{1}\framebox{1}\\ 
\framebox{2}
\end{array}$} 
\put(220,210){$\begin{array}{l}
\framebox{1}\framebox{3}\\ 
\framebox{2}
\end{array}$} 
\put(70,245){$\begin{array}{l}
\framebox{1}\framebox{1}\\
\framebox{3} \\
\end{array}$}
\put(170,262){$\begin{array}{l}
\framebox{1}\framebox{2}\\
\framebox{3}\\
\end{array}$}
\put(330,240){$\begin{array}{l}
\framebox{2}\framebox{2}\\
\framebox{3} \\
\end{array}$}
\put(210,373){$\begin{array}{l}
\framebox{1}\framebox{3}\\
\framebox{3} \\
\end{array}$}
\put(330,373){$\begin{array}{l}
\framebox{2}\framebox{3}\\
\framebox{3} \\
\end{array}$}
\put(152,120){\vector(-1,0){4}}
\put(272,360){\vector(-1,0){4}}
\put(90,182){\vector(0,-1){4}}
\put(330,302){\vector(0,-1){4}}
\put(274,184){\vector(-1,-1){1}}
\put(154,304){\vector(-1,-1){1}}
\put(205,310){\vector(-1,-4){1}}
\put(215,300){\vector(1,-4){1}}
\put(151,245){\vector(-4,-1){1}}
\put(280,235){\vector(-4,-1){1}}
\put(215,172){\vector(-1,-4){1}}
\put(264,304){\vector(-4,-3){1}}
\put(206,174){\vector(0,-1){1}}
\put(274,294){\vector(-3,-4){1}}
\put(150,171){\vector(-4,-3){1}}
\put(266,244.5){\vector(-4,-1){1}}
\put(200,174){1}
\put(268,365){1}
\put(146,306){1}
\put(335,298){1}
\put(221,298){1}
\put(196,307){1}
\put(260,310){1}
\put(280,290){-2}
\put(268,250){2}
\put(280,227){-1}
\put(150,250){1}
\put(83,180){1}
\put(150,112){1}
\put(280,182){-1}
\put(219,169){1}
\put(156,168){1}

\end{picture}}

${\mathbb S}^{\omega_1+\omega_2}({\mathbb C}^3)$
\end{center}

\

If we change our Weyl chamber, we can repeat this construction, defining first
anti semi standard tableaux as Young tableaux with entries strictly decreasing
in each column and decreasing in each row. Then we define an ordering on the set
of variables $\delta^{s}_{i_1,\dots,i_s}$, $i_1>i_2>\dots >i_s$ by putting:
$\delta^{(1)}>\delta^{(2)}>\dots$ and $\delta^{(p)}_{i_1,\dots,i_p}>
\delta^{(p)}_{j_1,\dots,j_p}$ if $i_p=j_p \dots i_{r+1}=j_{r+1}$ and $i_r>j_r$.\\

\

Let $T$ be an anti semi standard tableau. We can associate to $T$
a monomial:
$$\begin{array}{ccl}
\delta^T &=& \delta^{(c_1)}_{a^1_1\dots
a^1_{c_1}}\delta^{(c_2)}_{a^2_1 \dots a^2_{c_2}} \dots\\
&=&\pm \delta^{(c_1)}_{a^1_{c_1}\dots
a^1_1}\delta^{(c_2)}_{a^2_{c_2} \dots a^2_1} \dots \end{array}$$
and exchange the variables corresponding to columns with equal
height, then we get another Young tableau $T'$ such that $\delta^T
= \delta^{T'}$.\\
For instance:
$$T=\begin{array}{l}
\begin{array}{|c|c|}
\hline 4 & 2\\
\hline
\end{array}\\
\begin{array}{|c|c|}
\hline 3 & 1\\
\hline \end{array} \end{array},~~~~
\delta^{(2)}_{43}\delta^{(2)}_{21} =
\delta^{(2)}_{12}\delta^{(2)}_{34},~~~~ T'=\begin{array}{l}
\begin{array}{|c|c|}
\hline 1 & 3\\
\hline
\end{array}\\
\begin{array}{|c|c|}
\hline 2 & 4\\
\hline \end{array} \end{array} $$ or:
$$T=\begin{array}{l}
\begin{array}{|c|c|}
\hline 4 & 1\\
\hline
\end{array}\\
\begin{array}{|c|}
\hline 3 \\
\hline 2\\
\hline \end{array} \end{array},~~~~
\delta^{(3)}_{432}\delta^{(1)}_{1} =
-\delta^{(3)}_{234}\delta^{(1)}_{1},~~~~ T'=\begin{array}{l}
\begin{array}{|c|c|}
\hline 2 & 1\\
\hline
\end{array}\\
\begin{array}{|c|}
\hline 3\\
\hline 4\\
\hline\end{array} \end{array}. $$ Unfortunately, if $n>2$, $T'$ is
generally not semi standard $(T=\begin{array}{l}
\begin{array}{|c|c|}
\hline 3 & 1\\
\hline
\end{array}\\
\begin{array}{|c|}
\hline 2\\
\hline \end{array} \end{array},~~ T'=\begin{array}{l}
\begin{array}{|c|c|}
\hline 2 & 1\\
\hline
\end{array}\\
\begin{array}{|c|}
\hline 3\\
\hline \end{array} \end{array})$ thus our change of ordering on
the variables $\delta$ defines a new Groebner basis on the shape
algebra if $n>2$.\\
Now, the symmetry $\tau$ corresponds to the following operation on
tableaux since:
$$ \tau(\delta^{(s)}_{i_1,\dots,i_s})=\varepsilon^s_n \delta^{(s)}_{n+1-i_1,
\dots,n+1-i_s}$$
We can define $\tau$ directly on Young tableaux by replacing  each
entry $a^i_j$ of $T$ by $n+1-a^i_j$. The anti semi standard
tableaux are exactly the image by $\tau$ of the semi standard
ones.

\

\section{The reduced shape algebra : Algebraic presentation}

\

Let $V$ be a complex vector space with dimension $n$. From now one, we shall
study a quotient of the shape algebra  ${\mathbb S}^\bullet(V)$.

\begin{defi}

\

Let $R^+$ be the ideal in the shape algebra generated by $v^\lambda-1$:
\begin{align*}
R^+&=\langle v^\lambda-1=(e_1)^{a_1}(e_1\wedge e_2)^{a_2}\dots(e_1\wedge\dots
\wedge e_{n-1})^{a_n}-1,\quad\forall\lambda=\sum a_j\omega_j\rangle\\
&=\langle e_1-1,e_1\wedge e_2-1,\dots,e_1\wedge\dots\wedge e_{n-1}-1\rangle.
\end{align*}
We call \underbar{reduced shape algebra} and write ${\mathbb S}_{red}^\bullet(V)
$ the quotient $^{\hbox{${\mathbb S}^\bullet(V)$}}/_{\hbox{$R^+$}}$.

\end{defi}

\

This reduced shape algebra is no more a natural $\mathfrak{sl}(n)$
module but the ideal $R$ is invariant under the action of the
solvable group $HN^+$ consisting of upper triangula matrices in
$SL(n)$. Thus the quotient is a $HN^+$ module too. The action of
the Cartan group $H$ is still diagonal, let study the $N^+$ (or
$\mathfrak{n}^+$) action on ${\mathbb S}_{red}^\bullet(V)^+$.\\

\begin{prop}

\

Denote $\pi^+$ the canonical projection from ${\mathbb
S}^\bullet(V)$ to ${\mathbb S}_{red}^\bullet(V)^+$. Then

\begin{itemize}

\item i) The space of vectors $u\in{\mathbb S}_{red}^\bullet(V)^+$
such that $\mathfrak{n}^+u=0$ is ${\mathbb C}1$.

\item ii) ${\mathbb S}_{red}^\bullet(V)^+$ is an indecomposable
module.

\item iii) For any $\lambda$, the $\mathfrak{n}^+$ module
${\mathbb S}^\lambda(V )$ is equivalent to the submodule
$\pi^+\left({\mathbb S}^\lambda(V)\right)$ of ${\mathbb
S}_{red}^\bullet(V)^+$.

\item iv) For any $\lambda>\mu$, $\pi^+\left({\mathbb
S}^\mu(V)\right)$ is a submodule of $\pi^+\left({\mathbb
S}^\lambda(V)\right)$.
\end{itemize}

\end{prop}

\

\noindent
{\bf Proof}

\noindent {\it i)} We know (\cite{V} p. 317 for instance)
that, in each ${\mathbb S}^\lambda(V)$, the space of vectors $u$
such that $\mathfrak{n}^+u=0$ is exactly ${\mathbb C}v^\lambda$.
This gives {\it i)} in the quotient ${\mathbb
S}_{red}^\bullet(V)^+$.\\

\noindent {\it ii)} Let $u$ be a non zero vector in ${\mathbb
S}_{red}^\bullet(V)^+$, the $\mathfrak{n}^+$ module $W$ generated
by $u$ is finite dimensional since $u$ is a finite sum of image
through $\pi^+$ of weights vectors. The ${\mathfrak n}^+$ action
is locally nilpotent on ${\mathbb S}^\bullet(V)$, thus it is also
locally nilpotent on ${\mathbb S}_{red}^\bullet(V)^+$, as a
consequence $W$ contains a non trivial vector annhilated by
$\mathfrak{n}^+$. This vector is a multiple of 1. Thus any
$\mathfrak{n}^+$ submodule of ${\mathbb S}_{red}^\bullet(V)^+$
contains 1, ${\mathbb S}_{red}^\bullet(V)^+$ is an indecomposable 
$\mathfrak{n}^+$ module.\\

\noindent {\it iii)} Let $\pi^+_\lambda$ be the restriction of
$\pi^+$ to ${\mathbb S}^\lambda (V)$. It is a morphism of
$\mathfrak{n}^+$ modules. If its kernel is not vanishing, thanks
to Lie theorem, the $\mathfrak{n}^+$ module $Ker(\pi^+_\lambda)$
contains a non zero vector annihilated  by $\mathfrak{n}^+$, this
vector is a multiple of $v^\lambda$, but
$\pi^+(v^\lambda)=1\neq0$. Thus $\pi^+_\lambda$ is an
isomorphism of $\mathfrak{n}^+$ modules.\\

\noindent
{\it iv)} The relation $\lambda>\mu$ is equivalent to say there is $\nu$
dominant integral weight such that $\lambda=\mu+\nu$. In ${\mathbb S}^\bullet(V)
$, the multiplication by $v^\nu$ send ${\mathbb S}^\mu(V)$ into
${\mathbb S}^\lambda(V)$. In the quotient, this operation becomes the identity
mapping: $\pi^+(uv^\nu)=\pi^+(u)$ for any $u$ in  ${\mathbb S}^\mu(V)$.\\

\

Similarly, we define ${\mathbb S}_{red}^\bullet(V)^-$ as the
quotient of ${\mathbb S}^\bullet(V)$ by the ideal $R^-$ generated by
$\{e_{n}\wedge \dots \wedge e_{n+1-s}-1,~~s=1,\dots,n-1\}$. It is a
$HN^-$ module. If we denote $\pi^-$ the canonical morphism, we get
the same proposition with $''-''$ instead of $''+''$ everywhere.

\

\section{The reduced shape algebra, Geometrical presentation}

\

As above, we can write everything in term of the functions
$\delta^{(p)}_{i_1, \dots,i_p}$. If $R(\delta)^+$ is the ideal
generated by $\delta^{(p)}_{1,\dots,p} -1$, we get:
$${\mathbb S}_{red}^\bullet(V)^+\simeq^{\hbox{${\mathbb C}[SL(n,{\mathbb C})
]^{N^+}$}}/_{\hbox{$R(\delta)^+$}}=^{\hbox{${\mathbb
C}[\delta^{(p)}_{i_1,\dots,
i_p}]$}}/_{\hbox{$R(\delta)^++P(\delta)$}}.$$

Suppose now $f$ is a polynomial function, invariant with respect to the right
multiplication by $N^+$. Then $f$ is characterized by its restriction to the
dense open subset of $SL(n)$ of the matrices $g$ such that
$\delta^{(p)}_{1,\dots,p}(g)\neq0$ for all $p$. On this set, by the use of the
Gauss method, we can write:
$$g=\left[\begin{matrix}g'_{11}&0&0&\dots&0\\
g'_{21}&g'_{22}&0&\dots&0\\
g'_{31}&g'_{32}&g'_{33}&\dots&0\\
&&&\ddots&\\
g'_{n1}&g'_{n2}&g'_{n3}&\dots&g'_{nn}
\end{matrix}\right]\left[\begin{matrix}
1&a_{12}&a_{13}&\dots&a_{1n}\\
0&1&a_{23}&\dots&a_{2n}\\
0&0&1&\dots&a_{3n}\\
&&&\ddots&\\
0&0&0&\dots&1
\end{matrix}\right].$$
With, for all $k\geq j$:
$$g'_{jk}=\frac{\delta^{(k)}_{1,2,\dots,k-1,j}(g)}{\delta^{(k-1)}_{1,2,\dots,
k-1}(g)}.$$
By $N^+$ right invariance, we get
\begin{align*}
f(g)&=\frac{1}{\prod(\delta^{(j)}_{1,\dots,j}(g))^{b_j}}\Phi(
\delta^{(k)}_{1,2,\dots,k-1,j}(g),\delta^{(j)}_{1,\dots,j}(g))\\
&=\frac{1}{\prod(\delta^{(j)}_{1,\dots,j}(g))^{b_j}}\sum_{(c_1,\dots,c_{n-1})}
\Phi_{c_1,\dots,c_{n-1}}(\delta^{(k)}_{1,2,\dots,k-1,j}(g))\prod_j\left(
\delta^{(j)}_{1,\dots,j}(g)-1\right)^{c_j}\\
&=\frac{1}{\prod(\delta^{(j)}_{1,2,\dots,j}(g))^{b_j}}\sum_{(c_1,\dots,c_{n-1})}
F_{c_1,\dots,c_{n-1}}(g)\prod_j\left(\delta^{(j)}_{1,\dots,j}(g)-1\right)^{c_j}.
\end{align*}
By definition, the functions $\Phi_{c_1,\dots,c_{n-1}}$ and $F_{c_1,\dots,
c_{n-1}}$ are polynomial, $F_{c_1,\dots,c_{n-1}}$ is right invariant by $N^+$
and
$$F_{0,\dots,0}-f=\left(\prod(\delta^{(j)}_{1,\dots,j})^{b_j}-1\right)f-
\sum_{c_1+\dots+c_{n-1}>0}F_{c_1,\dots,c_{n-1}}\prod_j\left(\delta^{(j)}_{1,
\dots,j}-1\right)^{c_j}$$ belongs to $R(\delta)^+$. For any $g$ in
$N^-$ and any $k\geq j$, we have
$\delta^{(k)}_{1,\dots,k-1,j}(g)=g_{jk}$ and
$f(g)=\Phi_{0,\dots,0}(g_{jk})=F(g) $. The restriction of the
function $f$ to $N^-$ characterizes the function
$\Phi_{0,\dots,0}$ thus the value of $F_{0,\dots,0}$ and
$F_{0,\dots,0}$ and $f$ are in the same class modulo
$R(\delta)^+$. Conversely, any polynomial function $\Phi$ on $N^+$
defines a function $F$ in ${\mathbb C}[SL(n,{\mathbb C})]^{N^+} $.
The restriction mapping is an isomorphism of algebra between 
$\mathbb{S}^{\bullet}_{red}(V)^+$ and $\mathbb{C}[N^-]$.

\

\begin{rema}

\

In this presentation of ${\mathbb S}^{\bullet}_{red}(V)^+$, the
$N^+$ action on the elements of the reduced shape algebra is very
natural since it is just:
$$(g.f)(g') = f(^tgg'), ~~ g \in N^+, ~g' \in N^- ,~~ f \in \mathbb{C}[N^-]. $$
\end{rema}

\

But since $^{\hbox{${\mathbb C}[\delta^{(j)}_{i_1,\dots,i_j}]$}}
/_{\hbox{$R(\delta)^+$}}$ is simply ${\mathbb C}[\delta^{(j)}_{i_1,\dots,i_j}
(i_j>j)]$, we have also:
$$
{\mathbb S}_{red}^\bullet(V)^+\simeq^{\hbox{${\mathbb C}[\delta^{(j)}_{i_1,\dots
,i_j}(i_j>j)]$}}/_{\hbox{$P_{red}(\delta)^+$}}.
$$
Where $P_{red}(\delta)^+$ is the ideal generated by the Pl\"ucker relations but
where we replace the function $\delta^{(j)}_{1,\dots,j}$ by 1.

Especially, if $X_\alpha=E_{ij}$ $i<j$ then $X_\alpha$ acts on $\mathbb{C}
[\delta^{(p)}_{i_1,\dots,i_p}(i_p>p)]$ as the derivation:
\begin{align*}
X_\alpha f&=\frac{d}{ds}|_{s=0}f(exp~s^t X_\alpha.)\\&=
\pm\frac{\partial f}{\partial\delta^{(p)}_{(\{1,\dots,p\}\setminus\{i\})\cup\{
j\}}}+\sum\limits_{\begin{smallmatrix}\{i_1,\dots,i_p\}\cap\{i,j\}=\{j\}\\
\{i_1,\dots,i_p\}\setminus\{j\})\cup\{i\}\neq\{1,\dots,p\}\end{smallmatrix}}\pm
\delta^{(p)}_{(\{i_1,\dots,i_p\}\setminus\{j\})\cup\{i\}}\frac{\partial f}
{\partial\delta^{(p)}_{i_1,\dots,i_p}}.
\end{align*}

\

The same construction for ${\mathbb S}_{red}^\bullet(V)^-$ gives:
$$
R^-(\delta)=\theta\left(R^+(\delta)\right)
$$
is the ideal generated by the set $\{\delta^{(p)}_{n,\dots,(n-p+1)}-1\}$, 
${\mathbb S}_{red}^\bullet(V)^-$ is the quotient of ${\mathbb C}[SL(n)]^{N^+}$
(which is stabilized by $\theta$) by $R^-(\delta)$.  The Gauss formula allows to 
write:
$$
g=\left[ \begin{matrix}g'_{11}&&&&g'_{1n}\\
&&&.&\\&&.&&\\&.&&&\\g'_{nn}&&&&0\\
\end{matrix}\right]\left[\begin{matrix}
1&&&&a_{1n}\\
&.&&&\\
&&.&&\\
&&&.&\\
0&&&&1
\end{matrix}\right]
$$
if $\delta^{(p)}_{n,\dots,(n-p+1)}(g)\neq0$ for any $p$.

And any $f$ is modulo $R(\delta)^-$ characterized by its
restriction to:
$$
\left\{\left[ \begin{matrix}
g'_{ij}&&&&1\\
&&&.& \\
&&.&&\\
&.&&&\\
1&&&&0\\
\end{matrix} \right]\right\}=N^+\Omega\\
=\left\{\left[ \begin{matrix}
1&&&&g'_{ij}\\
&.&&& \\
&&.&&\\
&&&.&\\
0&&&&1\\
\end{matrix} \right]\left[ \begin{matrix}
0&&&&\varepsilon_n\\
&&&.& \\
&&.&&\\
&.&&&\\
\varepsilon_n&&&&0\\
\end{matrix}\right]\right\}
$$

Finally, if we put $f(n^+\Omega)=h(n^+)$, we get 
$\mathbb{S}^{\bullet}_{red}(V)^-\simeq \mathbb{C}[N^+]$ with the natural $N^-$ 
action:
$$
(g.h)(g_1)=h(^tgg_1).
$$

\begin{theo}

\

The reduced shape algebras are isomorphic to the algebra of
polynomial functions on $N^- (N^+)$, then:
$${\mathbb S}^\bullet_{red}(V)^+\simeq{\mathbb C}[N^-]={\mathbb C}[\mathfrak{n}^-]
$$
$${\mathbb S}^\bullet_{red}(V)^-\simeq{\mathbb C}[N^+]={\mathbb C}[\mathfrak{n}^+]
.$$
\end{theo}

\

The last assertions of the theorem comes from the observation that
the exponential mapping from the Lie algebra
$\mathfrak{n}^-$($\mathfrak{n}^+$) onto the Lie group $N^-$($N^+$)
is a polynomial bijection with inverse polynomial too.

\

\section{The reduced shape algebra: Combinatorial presentation}

\subsection{Super and quasi standard Young tableaux}

\

In order to describe the restricted shape algebra and the
restricted Pl\"ucker relations, we have to perform the quotient of
the preceding construction by the ideal generated by
$\{\delta^{(s)}_{12\dots s}-1\}$. On the Young tableaux this
operation can be viewed as an 'extraction' of trivial columns.\\

A column whose height is $c$ in a tableau is \underbar{trivial} if
its entries are $1,2,\dots,c$, a Young tableau $T$ is trivial if
each column of $T$ is trivial. Now let $T$ be a Young tableau
(semi standard or not), we define the extraction of trivial
columns in $T$ in the following manner:

Denote $a_{ij}$ the entries of $T$ ($a_{ij}$ is in the row $i$ and
the column $j$, for any $j$, $a_{ij}<a_{(i+1)j}$ and the heights
$c_1,\dots,c_t$ of $T$ are decreasing). We say that the tableau
$T$ is \underbar{reducible} if
\begin{itemize}
\item there is a column $j$ whose the $s$ top entries are $1,2,\dots,s$
($a_{i,j}=i$ for $1\leq i\leq s$),

\item on the right of the column $j$, there is a column $j'$ with height $s$
in $T$ (there is $j'\geq j$ such that $c_{j'}=s$),

\item for any $k>j$, if $c_{k-1}>s$ and $c_k\geq s$,
$a_{s+1,k-1}>a_{s,k}$.
\end{itemize}

Let $T$ be a reducible Young tableau, let $j$ the smallest index
and $s$ the largest index for which the above conditions hold. Let
us suppress the trivial top part of the column $j$ and shift to
the left the right parts of the $s$ first rows ({\sl i.e.} we
shift to the left every $a_{ik}$ with $1\leq i\leq s$ and $j<k$),
then we get a Young tableau $R_1$: the entries of $R_1$ are
$b_{k\ell}$ with
$$b_{k\ell}=\left\{\begin{matrix}a_{k(\ell+1)}&\hbox{ if }&1\leq k\leq s&\hbox{
and }&j\leq\ell\leq t-1\cr
&&&&\cr
a_{k\ell}&\hbox{ if }&s<k&\hbox{ or }&\ell<j.
\end{matrix}\right.$$
If the number of column of $T$ was $t$, then $R_1$ has $t-1$ column, more
precisely if the heights of the columns of $T$ were: $(c_1,\dots,c_t)$ and the
columns of heights $s$ had the number $j',\dots,j''$, then the heights of the
columns of $R_1$ are $(c'_1,\dots,c'_{t-1})$ with
$$c'_k=\left\{\begin{matrix}c'_k=c_k&\hbox{ if }&1\leq k<j''\cr
&&\cr c'_k=c_{k+1}&\hbox{ if }&j''\leq k\leq
t-1.\end{matrix}\right.$$ Simultaneously, we define $L_1$ as the
Young tableau with only one trivial
column with entries $1,\dots,s$.\\

Now if $R_1$ is reducible, we repeat the above operation,
extracting a second trivial column from $R_1$, getting two Young
tableaux a trivial one with two columns $L_2$ and a Young tableau
$R_2$ with $t-2$ columns.

Repeating this construction, after $m$ steps, we get a trivial Young tableau
$L_m$ with $m$ column and a Young tableau $R_m$ with $t-m$ columns.

This construction stops when the Young tableau $R_m$ is not reducible we say
$R_m$ is irreducible and call $R_m$ the residue of $T$.\\

\

\

\begin{defi} {\rm (Super, left and right Young tableaux)}

A \underbar{super Young tableau} is a pair $S=(L,R)$ of two Young tableaux, the
left one $L$ is a trivial Young tableau, the right one, $R$ is an irreducible
Young tableau. $L$ or $R$ can be the empty tableau without any column.
\end{defi}

\

Our construction defines a mapping $f$ (the extraction mapping) from the set
$\mathcal{Y}$ of Young tableaux into the set $\mathcal{SY}$ of super young
tableaux
$$f(T)=S=(L,R).$$

If $\lambda$ is the sequence of the heights of the column of $T$:
$\lambda=(c_1\leq \dots\leq c_t)$ and $\mu=(c'_1\leq\dots\leq
c'_\ell)$ and $\nu=(c''_1\leq \dots\leq c''_r)$ the corresponding
sequence for $L$ and $R$ (one of these sequences can be empty),
then $\mu$ and $\nu$ are two disjoint subsequences of $\lambda$
and $\lambda$ is the 'union' of $\mu$ and $\nu$: we refind the
sequence $\lambda$ by putting together the elements of $\mu$ and
$\nu$ and
ordering them in a decreasing sequence.\\

Starting with an irreducible Young tableau $R$, we can insert to
it any family of trivial columns, say $\{D_1,\dots,D_\ell\}$,
getting a new tableau $T$. We insert these column in the following
way: if the height of $D_i$ is $d_i$, we insert $D_1,\dots,D_i$
such that any column of $T$, after $D_i$ has height strictly less
then $d_i$, the columns of $T$ before $D_i$ are the columns of $R$
with length at least $d_i$, with their ordering and the column
$D_j$ ($j<i$). Then $T$ is a Young tableau. Of course, if
$\ell>0$, $T$ is reducible.

\

Let us try to extract a trivial column from $T$. Among the new
column, the first one is $D_1$ with height $d_1$. In $T$ this
column is the column $p$. Suppose the first trivial column
extracted from $T$ is the $s$ top elements of the column $j$, with
$j<p$. Since $R$ is irreducible, there is a $k>j$ such that
$c_{k-1}^R>s$, $c_k^R\geq s$ and $a^R_{s+1,k-1}\leq a^R_{s,k}$ (we
denote $c_k^R$ the height of the column $k$ and $a_{i,j}^R$ the
$i,j$-entry in $R$). We choose the smallest such $k$. Since we can
now extract the trivial column from $T$, there is, in $T$, at
least one new column, say $D$ between the two columns $k-1$, $k$
in $R$, which are now columns $k_1$, $k_2$ in $T$. We choose for
$D$ the last one: $D$ is the column $k_2-1$ in $T$. The height of
$D$ is $c^T_{k_2-1}=d>c^T_{k_2}= c^R_k\geq s$ and we get:
$$a^T_{s+1,k_2-1}=s+1\leq a^R_{s+1,k-1}\leq a^R_{s,k}=a^T_{s,k_2}.$$
We cannot extract the trivial column consisting of the $s$ top elements of the
column $j$, with $j<p$. Of course, we can extract all the column $p$ of $T$.
Thus, in the computing of $f(T)$, the first step is just to eliminate the column
$D_1$ from $T$, repeating this construction, we get $f(T)=(L,R)$ where $L$ is
the trivial tableau $(D_1,\dots,D_\ell)$. We proved that $f$ is a surjective
mapping by defining a mapping $h$ from $\mathcal{SY}$ to $\mathcal{Y}$ such that
$f\circ h(L,R)=(L,R)$.\\

\begin{defi} ({\rm Quasi standard tableaux)}

A super Young tableau $S=(L,R)$ is said \underbar{quasi-standard} if its right
tableau $R$ is semistandard.

A Young tableau $T$ is said quasi-standard if it is irreducible and
semistandard.
\end{defi}

Let us denote by $\mathcal{QSY}$ (resp. $\mathcal{QY}$) the set of quasi
standard super Young tableaux (resp. quasi standard Young tableaux). Denote
$\mathcal{SEM}$ the set of semistandard Young tableaux.\\

\begin{lem} ({\rm $f$ is a bijection from $\mathcal{SEM}$ onto $\mathcal{QSY}$)}

The mapping $f$, when restricted to $\mathcal{SEM}$ is a
one-to-one onto mapping from $\mathcal{SEM}$ onto $\mathcal{QSY}$.
\end{lem}

\

\noindent
{\bf Proof}

First it is clear that if $T$ is semistandard, then each tableau in the sequence
 $R_1,\dots,R_m$ defined above is still semistandard, then $f$ is a map from
$\mathcal{SEM}$ to $\mathcal{QSY}$.\\

Now let $S=(L,R)$ be an element of $\mathcal{QSY}$. Denote the
rows of $L$ by $(L'_1,\dots,L'_u)$, their lengths being
$\ell'_1,\dots,\ell'_u$. Similarly, denote $(L''_1,\dots,L''_v)$
the rows of $R$, their lengths being $\ell''_1,\dots, \ell''_v$.We
define the new tableau $T=g(S)$ as the tableau with the row $i$
contains (from left to right) $\ell'_i$ entries $i$, then the
$\ell_i''$ entries of the row $i$ of $R$. In fact, $T$ is a Young
tableau since if $a^T_{i,j}$ is an entry of $T$, it is either $i$
or an entry of $R$ ($a^T_{i,j}=a^R_{i,j- \ell'_i}$ if $a^R_{r,s}$
are the entries of $R$). In any case, $a^T_{i,j}\geq i$.

If $a^T_{i,j}=i$, then $a^T_{i,j}=i<i+1\leq a^T_{i+1,j}$. If $a^T_{i,j}=a^R_{i,
j-\ell'_i}$, since $\ell'_i\geq \ell'_{i+1}$, $a^T_{i+1,j}=a^R_{i+1,j-
\ell'_{i+1}}$ and $a^T_{i,j}=a^R_{i,j-\ell'_i}<a^R_{i+1,j-\ell'_i}\leq a^R_{i+1,
j-\ell'_{i+1}}=a^T_{i+1,j}$.

$T$ is semistandard: by construction each row in $T$ is a increasing sequence of
entries. $g$ is a map from $\mathcal{QSY}$ to $\mathcal{SEM}$.\\

The map $g$ is the inverse mapping of $f$. Indeed if $T$ is semistandard, if a
column $C$ of $T$ begins by a trivial part, then all the columns before $C$
begin with the same trivial part and suppressing the top of the first column or
the top of $C$ is the same operation, thus to construct the sequence $R_1,
\dots, R_m$, we just have to consider the first column at each step.

Starting with $T=g(S)$, we can extract at each step a trivial column having the
height of the corresponding column of $L$, but no more, since $R$ is
irreducible. Thus $f\circ g(S)=S$, for any $S\in\mathcal{QSY}$.

Conversely, starting with a semistandard $T$, we build first
$f(T)=(L,R)$ and by construction the rows of $L$ are the left part
of the rows of $T$, thus $g\circ f(T)=T$.\\

\

\subsection{Quasi standard Young tableaux and Groebner basis}

\

In this section we shall repeat the construction of section 3 but for the ideal
$R(\delta)^+$ and the quasi standard Young tableaux.\\

First, we choose the following elimination order on the variables $\delta$:
defining the degree $deg(\delta^{(s)}_{i_1\dots i_s})$ as 1 if $i_s>s$
($\delta^{(s)}_{i_1\dots i_s}$ is not trivial) and 0 if $i_s=s$
($\delta^{(s)}_{i_1\dots i_s}$ is trivial), the degree of $\delta^T$ is the sum
of degree of each variables and $T>T'$ if and only if:
$$\left\{\begin{matrix}
&deg(\delta^T)>deg(\delta^{T'})&&\cr
\hbox{or}&&&\cr
&deg(\delta^T)=deg(\delta^{T'})&\hbox{and}&T>T'~\hbox{for the preceding
ordering}.\end{matrix}\right.$$

Now we look for the leading terms of elements of $R(\delta)^+$, for
this ordering. We saw that the leading terms of elements of
$P(\delta)$ for the preceding
ordering were non semistandard monomials.\\

Let $T$ be a  non quasi standard tableau.\\

\noindent
\underbar{\bf Case 1}: $T$ is non semi standard.\\

Then $T$ contains a non semistandard tableau with two columns
$T^0$: $\delta^T= \delta^U\delta^{T^0}$. For $T^0$, we saw there
is a Pl\"ucker relation $P_{T^0}$ in $P(\delta)$ whose leading
term for the ordering of section 3 was $T^0$.\\

\noindent
\underbar{Case 1.1}: $T^0$ contains a trivial column
$C_i$, since $T^0$ is non semistandard, it is its second column.
$\delta^{T^0}=\delta^{(s)}_{1,\dots,s}
\delta^{(c)}_{a_1,\dots,a_c}$. But $\delta^{(s)}_{1,\dots,s}$ is
the leading term of the element $V_s=\delta^{(s)}_{1,\dots,s}-1$
in $R(\delta)^+$. $\delta^{T^0}$ is the leading term of
$\delta^U\delta^{(c)}_{a_1,\dots,a_c}V_s$ which is in $R(\delta)^+$.\\

\noindent
\underbar{Case 1.2}: $T^0$ does not contain any trivial column. $\delta^{T^0}=
\delta^{(s)}_{b_1,\dots,b_s}\delta^{(c)}_{a_1,\dots,a_c}$ with $c\geq s$, there
is $j$ such that $a_j> b_j$, we choose the largest such $j$, due to our
conventions of writing, if $c=s$ then $j<s$ and $a_c>c$, $b_s>s$.

Thus the relation $P_{T^0}$ has the following form:
$$\begin{array}{cl}
P_{T^0}&=\delta^{T^0}-\sum\limits_{\begin{smallmatrix}
A\subset\{a_1,\dots,a_c\}\cr \# A=j\end{smallmatrix}} \pm
\delta^{(s)}_{A\cup \{b_{j+1},\dots,b_s\}}
\delta^{(c)}_{(\{a_1,\dots,a_c\}\setminus
A)\cup\{b_1,\dots,b_j\}}\\
&=\delta^{T^0}-\sum\limits_{\begin{smallmatrix} S<T^0\cr S~{\rm
semi~standard}
\end{smallmatrix}} \pm \delta^{(S)}.\end{array}$$
If a tableau $S$ in this relation contains a trivial column, {\sl i.e}
$S=C_1C_2$ with $C_1$ trivial, we replace $S$ by $C_2$ since
$$
\delta^S-\delta^{C_2}=V_s.\delta^{C_2}.
$$
Repeating this operation, we get an element
$$
P'_{T^0}=\delta^{T^0}-\sum\limits_{\begin{smallmatrix} S<T^0\cr
S~{\rm quasi~standard}
\end{smallmatrix}} \pm \delta^{(S)}.
$$\\

\noindent
\underbar{\bf Case 2}: $T$ is semi standard.\\

If $T$ has only one column, this column is trivial $T$ is the
leading term of some $P_T = \delta^T - 1$ in $R(\delta)^+$.

Since $T$ is semi standard the construction of the super Young tableau $f(T)$
begins with the extraction of the top $s$ elements $1,\dots,s$ of the first
column of $T$. Let us look to the two first columns of $T$, $C_1^T$ and $C_2^T$.
By hypothesis, $\delta^{C_1^T}=\delta^{(c_1)}_{1,\dots,s,a_{s+1},\dots,a_{c_1}}
$, $\delta^{C_2^T}=\delta^{(c_2)}_{b_1,\dots,b_s,b_{s+1},\dots,b_{c_2}}$ and
$b_s<a_{s+1}$.

Let us define $\partial T$ as the tableau with the following first columns 
$C_1^{\partial T}$ and $C_2^{\partial T}$:
$$
\delta^{C_1^{\partial T}}=\delta^{(c_1)}_{b_1,\dots,b_s,a_{s+1},\dots,a_{c_1}}
,~~~\delta^{C_2^{\partial T}}=\delta^{(c_2)}_{1,\dots,s,b_{s+1},\dots,b_{c_2}},
$$
the other columns of $\partial T$ being $C_i^{\partial T}=C_i^T$
($i\geq 3$). Let us write the Pl\"ucker relation corresponding to
these two columns and $s$:
\begin{align*}
&\delta^T-\delta^{\partial T}-\\
&\sum_{\begin{smallmatrix}A\subset\{1,\dots,s,a_{s+1},\dots, a_{c_1}\}\\
A\neq\{1,\dots,s\}\\
\# A=s\end{smallmatrix}}\pm
\prod_{i\geq3}\delta^{C_i^T}\delta^{(c_2)}_{A\cup\{b_{s+1}
,\dots,b_{c_2}\}}\delta^{(c_1)}_{(\{1,\dots,s,a_{s+1},\dots,a_{c_1}\}\setminus
A )\cup\{b_1,\dots, b_s\}}\\
&=\delta^T-\delta^{\partial T}-\sum_A \pm \delta^{T_A}.
\end{align*}

Each term $\delta^{T_A}$ in the sum has a second
column containing $a_i$ with $i>s$, thus $a_i\geq a_{s+1}>b_s$ and
$\delta^{(c_2)}_{A\cup\{b_{s+1},\dots, b_{c_2}\}}<\delta^{C_2^T}$,
$\delta^{T_A}<\delta^T$.

If $c_2=s$, $deg(\delta^{\partial T})<deg(\delta^T)$, $\delta^T$
is the leading term of an element in $R(\delta)^+$. If $c_2>s$, we
repeat this construction for $\partial T$, forgotting its first
column. We get the following element of $R(\delta)^+$:
\begin{align*}
&\delta^{\partial T}-\delta^{\partial^2 T}-\\
&\sum_{\begin{smallmatrix}
B\subset\{1,\dots,s,b_{s+1},\dots, b_{c_2}\}\\
B\neq\{1,\dots,s\}\\
\# B=s\end{smallmatrix}} \pm
\prod_{i\geq4}\delta^{C_i^T}\delta^{(c_3)}_{B\cup
\{c_{s+1},\dots,c_{c_3}\}}\delta^{(c_2)}_{(\{1,\dots,s,b_{s+1},\dots,b_{c_2}\}
\setminus B)\cup\{c_1,\dots, c_s\}}\delta^{C_1^{\partial T}}\\
&= \delta^{\partial T} - \delta^{\partial^2 T}- \sum_B
\delta^{T_B}.
\end{align*}
Each term $\delta^{T_B}$ in the sum has a third column containing $b_i$ with
$i>s$, thus $b_i\geq b_{s+1}>c_s$ and $\delta^{(c_3)}_{B\cup\{c_{s+1},\dots,
c_{c_3}\}}<\delta^{C_3^T}$, $\delta^{T_B}<\delta^T$.

Repeating this operation we finally get an element in $R(\delta)^+$ of the form:
$$
\delta^T - \delta^{\partial^k T}-\sum_j \delta^{T_j}$$ with
$\delta^{T_j}<\delta^T$ for all $j$, the column $k+1$ of
$\partial^k T$ is trivial, $deg(\delta^{\partial^k
T})<deg(\delta^T)$ and $\delta^T$ is the
leading term of an element of $R(\delta)^+$.\\

\begin{rema}

\

The tableau $\partial^kT$ considered here is (perhaps up a
reordering of the columns with height $s$ ) the tableau $h(C,R_1)$
if $C$ is the first trivial column: $$ \delta^C = \delta^{(s)}_{1
\dots s}$$ and $R_1$ the first step in the process of trivial
columns extraction from $T$.\end{rema}
We got an element of
$R(\delta)^+$: 
$$
\delta^T - \delta^{R_1} - \sum \pm \delta^{T_j}
~~~ (R_1 < T, T_j < T)
$$
If $R_1$ is quasi standard, we stop the process. If it is not the
case, we continue the extraction, getting new tableaux $T'_k < R_1
< T$. Finally we get: $ T = g(L, R) ~ {\rm with} ~L \neq \emptyset
~{\rm and}~ \delta^T - \delta^R - \sum\limits_{T_k < T} a_k
\delta^T_k $ belongs to $R(\delta)^+$. $R$ is quasi standard $R < T,
~T_k < T$. We repeat this operation for each non quasi standard
$T_k$; getting an element $P_T = \delta^T - \delta^R - \sum a_k
\delta^{T_k}$ with $T_k < T$, $T_k$ quasi standard, $P(T)$ in
$R(\delta)^+$.

We proved that each non quasi standard Young tableau is the
leading term of an explicit element $P_T$ of $R(\delta)^+$. Let us
now prove that any quasi standard Young tableau is not a leading
term of an element in $R(\delta)^+$.\\

Let $\lambda$ be a highest weight for $\frak{sl}(n)$ and
$V^\lambda$ the corresponding simple module. We saw that
$V^\lambda$ is naturally a sub-module of ${\mathbb
S}^\bullet_{red}(V)$. More precisely, $V^\lambda$ is the space
spanned by the classes modulo $R(\delta)^+$ of the monomials
$\delta^T$ for all Young tableau $T$ of shape $\lambda$. A basis
for $V^\lambda$ is given by the classes of the monomials
$\delta^T$ for $T$ semi standard with shape $\lambda$ in the
quotient ${\mathbb C}[\delta]/R(\delta)^+$. Let us consider the
sub-space $W^\lambda$ of $V^\lambda$ spanned by the quasi standard
and semi standard Young tableau of shape $\lambda$.\\
A basis of $V^{\lambda}$ is given by the classes of $\delta^T$,
$T$ semi standard with shape $\lambda$ modulo $R(\delta)^+$. Either
$T$ is quasi standard and $(\delta^T)$ is a basis of $W^{\lambda}$
or $T=g(L,R)$, we saw $\delta^T-\delta^R=\sum a_k\delta^{T_k}$
modulo $R(\delta)^+$ with $T_k$ quasi standard $T_k<T$ and $R$ is
quasi standard with shape $\mu<\lambda$. This proves that
$V^\lambda$ is a subspace of
$\sum\limits_{\mu\leq\lambda}W^\mu$.But since $g$ is injective,
$dim V^\lambda = \sum\limits_{\mu\leq\lambda}dim(W^\mu)$ thus
$V^\lambda=\displaystyle{\bigoplus_{\mu\leq\lambda}W^\mu}$.

Let now $T$ be a quasi standard Young tableau of shape $\lambda$.
Suppose $\delta^T$ is the leading term of an element
$T+\sum_ka_k\delta^{T_k}$ in $R(\delta)^+$, then using the first
part of the proof, we can replace each $\delta^{T_k}$ with a non
quasi standard $T_k$,  by a linear combination of $\delta^{T_j}$
with quasi standard $T_j$ modulo $R(\delta)^+$. Finally we get an
element in $R(\delta)^+$ of the form $T+\sum_ja_j\delta^{T_j}$ with
any $T_j$ quasi standard and strictly smaller than $T$. The shape
$\mu_j$ of $T_j$ is thus smaller than $\lambda$. But this is
impossible since the sum $\sum\limits_{\mu\leq\lambda}W^\mu$ is
direct.

Finally, as in section 3, for each non quasi standard Young
tableau, we got an element in $R(\delta)^+$ of the form:
$$P^{red}_T=\delta^T-\sum_ja_j\delta^{T_j}$$
with $\delta^{T_j}$ strictly smaller than $\delta^T$ and quasi
standard.\\

Let $T$ be a non quasi standard tableau with shape
$\lambda$. We shall say that $T$ is minimal if it does not contain
any non quasi standard tableau with shape $\mu < \lambda$. For
instance a semi standard non quasi standard tableau with one
column or with 2 columns without trivial column are minimal.

If $n\leq 3$ there are no other semi standard, minimal, non quasi
standard tableaux, but if $ n \geq 4$ there is semi standard,
minimal, non quasi standard tableau with 3 columns for instance:

\begin{displaymath}
\begin{array}{c}
\begin{array}{|c|c|c|}
\hline
1&2&3\\
\hline
\end{array}\\
\begin{array}{|c|c|}
3&4\\
\hline
\end{array}\hfill\\
\end{array}
\end{displaymath}

\

\begin{theo} {\rm The Groebner basis}

\

The set
\begin{align*}
G&=\{P^T_{red},~T~{\it semi~standard~minimal~non~quasi~standard~or}\\
&~T~ {\it non~semi~standard~with~2~columns,~without~any~trivial~column}\}
\end{align*}
is the reduced Groebner basis of $R(\delta)^+$ for our ordering.
\end{theo}

\

\noindent
{\bf Proof:}

We saw that
$$
<LT(R(\delta)^+)>=\{\delta^T, T {~\rm non~ quasi~standard}\}.
$$

If $T$ contains a trivial column $C$, $\delta^T$ is divisible by
$\delta^C$ and $C$ is minimal semi standard non quasi standard.\\

If $T$ is non semi standard, it contains a non semi standard tableau $S$ with 2 
columns, without any trivial column.\\

If $T$ does not contain any trivial column and is semi standard
then by definition it contains a minimal non quasi standard
tableau $S$ but $S$ is by construction semi standard.\\

Thus
$$
<LT(G)>=<LT(R(\delta)^+)>.
$$
Now each monomial in any $P^T_{red}$ of $G$ which are not the leading term, are 
$a_{T'}\delta^{T'}$ with $T'$ quasi standard.

But if $S\subset T'$, then $S$ is also quasi standard. Indeed, $S$ is semi 
standard, suppose $S$ non quasi standard then $S$ contains a first column
$$
C_1=(1, 2, \dots, s, a_{s+1}, \dots, a_{C_1}),
$$
other column
$$
C_i=(b_1, b_2, \dots, b_s, b_{s+1}, \dots)
$$
a last column
$$
C_t=(c_1, c_2, \dots, c_{c_t})~~\hbox{with}~~t\leq s.
$$
We can extract $(1, 2, \dots, s)$ from $S$.

Now, we can refind $T$ from $S$ by adding some columns before $C_1$, between 
columns of $S$ or after $C_t$. But $T$ is semi standard. By considering each 
case for these new columns, we directly see that the top $(1, 2, \dots, s)$ of 
columns $C_1$ can still be extracted from $T$ which is impossible since $T$ is
quasi standard.\\

Thus any monomial of $P^T_{red}$ is not divisible by the leading
term of another $P^T_{red}$.\\

This means that $G$ is the reduced Groebner basis of $R(\delta)^+$ for our 
ordering.

\

The same result holds with the anti standard tableau, image by $\tau$ of the 
quasi standard tableaux.

\

The anti quasi standard tableaux can be defined exactly as the
quasi standard tableaux by extracting ''trivial'' top of columns like:
\begin{displaymath}
\begin{array}{|c|}
\hline
n\\
\hline
n-1\\
\hline
\vdots\\
\hline
n-s\\
\hline
\end{array}
\end{displaymath}

They are still the image by $\tau$ of the quasi standard tableaux.

\

\begin{rema}

\

In fact, if $n\leq3$, the quasi satndard Grobner basis is invariant under the 
action of $\theta$. Similarly, with the symmetry $\tau$, if we identify 
$\tau(T)$ with $\pm T'$ with $T'$ the Young tableau such that $\deltaç{\tau(T)}=
\deltaç{T'}$, then $T$ quasi standard  implies $T'$ quasi standard. In the study
of ${\mathfrak sl}(4)$ below, we shall see this no more true for $n>3$.
\end{rema}

\

Let us now picture the adjoint representation of $\frak{sl}(3)$ in 
$\mathbb{S}^+_{red}$ equipped with its Groebner basis:

\begin{center}

{\tiny
\begin{picture}(410,600)(0,0)
\put(90,120){\vector(1,0){340}} \put(90,120){\vector(0,1){320}}
\put(100,430){$(2,3)$} \put(410,110){$(1,2)$}
\path(210,120)(330,240)
\path(330,240)(330,360)
\path(330,360)(210,360)
\path(210,360)(90,240)
\dottedline[.](330,360)(200,250)
\dottedline[.](330,360)(220,230)
\dottedline[.](210,360)(200,250)
\dottedline[.](210,360)(220,230)
\dottedline[.](330,240)(200,250)
\dottedline[.](330,240)(220,230)
\dottedline[.](90,240)(200,250)
\dottedline[.](90,120)(200,250)
\dottedline[.](90,120)(220,230)
\dottedline[.](210,120)(220,230)
\put(90,120){\circle{4}} \put(90,240){\circle{4}}
\put(90,360){\circle{4}} \put(210,120){\circle{4}}
\put(330,120){\circle{4}} \put(200,250){\circle{4}}
\put(220,230){\circle{4}} \put(330,240){\circle{4}}
\put(330,360){\circle{4}} \put(210,360){\circle{4}} \put(202,110){
\framebox{2}} \put(220,220){ \framebox{3}} \put(85,110){0}
\put(318,110){ \framebox{2}\framebox{2}}
\put(70,240){$\begin{array}{l}
\framebox{1}\\
\framebox{3} \\
\end{array}$}
\put(60,360){$\begin{array}{l}
\framebox{1}\framebox{1}\\
\framebox{3}\framebox{3} \\
\end{array}$}
\put(180,262){$\begin{array}{l}
\framebox{2}\\
\framebox{3} \\
\end{array}$}
\put(330,240){$\begin{array}{l}
\framebox{2}\framebox{2}\\
\framebox{3} \\
\end{array}$}
\put(210,373){$\begin{array}{l}
\framebox{1}\framebox{3}\\
\framebox{3} \\
\end{array}$}
\put(330,373){$\begin{array}{l}
\framebox{2}\framebox{3}\\
\framebox{3} \\
\end{array}$}
\put(152,120){\vector(-1,0){4}} \put(272,360){\vector(-1,0){4}}
\put(90,182){\vector(0,-1){4}} \put(330,302){\vector(0,-1){4}}
\put(274,184){\vector(-1,-1){1}} \put(154,304){\vector(-1,-1){1}}
\put(205,310){\vector(-1,-4){1}} \put(215,300){\vector(1,-4){1}}
\put(151,245){\vector(-4,-1){1}} \put(280,235){\vector(-4,-1){1}}
\put(215,172){\vector(-1,-4){1}}
\put(264,304){\vector(-4,-3){1}} \put(145,186){\vector(-3,-2){1}}
\put(274,294){\vector(-3,-4){1}} \put(150,171){\vector(-4,-3){1}}
\put(266,244.5){\vector(-4,-1){1}} \put(268,365){1}
\put(146,306){1} \put(335,298){1} \put(221,298){2}
\put(196,307){1} \put(260,310){1} \put(280,290){-1}
\put(268,250){2} \put(280,227){1} \put(150,250){1} \put(83,180){1}
\put(150,112){1} \put(280,182){-1} \put(219,169){1}
\put(137,188){-1} \put(156,168){1}
\end{picture}}

\end{center}

and the same representation in $\mathbb{S}^-_{red}$ equipped with its Groebner
basis is:
\begin{center}
{\tiny
\begin{picture}(410,600)(0,0)
\put(330,500){\vector(-1,0){320}}
\put(330,500){\vector(0,-1){320}}
\put(340,180){$(2,3)$}
\put(20,510){$(1,2)$}

\path(210,500)(90,380)
\path(330,380)(210,260)
\path(210,260)(90,260)
\path(90,380)(90,260)

\dottedline[.](330,500)(200,390)
\dottedline[.](330,500)(220,370)
\dottedline[.](210,500)(220,370)
\dottedline[.](330,380)(200,390)
\dottedline[.](210,260)(200,390)
\dottedline[.](210,260)(220,370)
\dottedline[.](90,260)(200,390)
\dottedline[.](90,260)(220,370)
\dottedline[.](90,380)(220,370)
\dottedline[.](90,380)(200,390)

\put(330,500){\circle{4}} \put(210,500){\circle{4}}
\put(330,380){\circle{4}} \put(200,390){\circle{4}}
\put(220,370){\circle{4}} \put(90,380){\circle{4}}
\put(90,260){\circle{4}} \put(210,260){\circle{4}}
\put(90,500){\circle{4}} \put(330,260){\circle{4}}

\put(335,260){ \framebox{2}\framebox{2}} \put(335,380){
\framebox{2}}  \put(335,505){0} \put(200,515){$\begin{array}{l}
\framebox{3}\\
\framebox{1} \\
\end{array}$}
\put(75,515){$\begin{array}{l}
\framebox{3}\framebox{3}\\
\framebox{1}\framebox{1} \\
\end{array}$}
\put(60,370){$\begin{array}{l}
\framebox{3}\framebox{1}\\
\framebox{1} \\
\end{array}$}
\put(180,400){ \framebox{1}} \put(220,355){$\begin{array}{l}
\framebox{2}\\
\framebox{1} \\
\end{array}$}
\put(200,245){$\begin{array}{l}
\framebox{2}\framebox{2}\\
\framebox{1} \\
\end{array}$}
\put(60,245){$\begin{array}{l}
\framebox{2}\framebox{1}\\
\framebox{1} \\
\end{array}$}

\put(148,260){\vector(1,0){4}} \put(268,500){\vector(1,0){4}}
\put(90,322){\vector(0,1){4}} \put(330,440){\vector(0,1){4}}
\put(270,320){\vector(1,1){1}} \put(150,440){\vector(1,1){1}}
\put(262,442){\vector(3,4){1}} \put(264,385){\vector(4,-1){1}}
\put(144,385){\vector(4,1){1}} \put(215,314){\vector(0,4){1}}
\put(154,375){\vector(4,-1){1}}

\put(215,440){\vector(0,1){1}} \put(205,327){\vector(0,1){1}}
\put(274,434){\vector(3,4){1}} \put(144,324){\vector(3,4){1}}
\put(154,314){\vector(4,3){1}}

\put(149,253){1} \put(85,320){1} \put(270,504){1} \put(219,316){2}
\put(333,440){1} \put(146,442){1} \put(272,316){-1}
\put(144,388){2} \put(262,378){1} \put(152,370){1}
\put(211,440){1} \put(160,312){1} \put(279,434){-1}
\put(258,444){1} \put(136,324){-1} \put(200,325){1}

\end{picture}}

\end{center}

We resume our construction by the two following diagrams:

$$\begin{array}{c}
A^\bullet(V)=\mathbb{C}[V\oplus(V\wedge V)\oplus\dots\oplus(V\wedge\dots\wedge 
V)]\\
\vspace{.1cm}\\
\downarrow\\
\vspace{.1cm}\\
\tau\curvearrowleft\\
\mathbb{S}^{\bullet}(V)=^{\mathbb{C}[V\oplus \dots \oplus\wedge^{n-1}V]}/_P\\
\vspace{.1cm}\\
\pi^-\swarrow~~~~~~~~~~\searrow\pi^+\\
\vspace{.1cm}\\
^{\mathbb{C}[V\oplus \dots \oplus\wedge^{n-1}V]}/_{(P+R^-)=\mathbb{S}^-_{red}}
~~~~\longleftrightarrow~~~~
\mathbb{S}^+_{red}=^{\mathbb{C}[V\oplus \dots \oplus\wedge^{n-1}V]}/_{(P+R^+)}
\end{array}$$
and :
$$\begin{array}{c}
\theta\curvearrowleft\\
\begin{array}{rcl}
^{\mathbb{C}[\delta^{(s)}_{i_1, \dots,i_s}]_{i_1>\dots>i_s}}/_{P(\delta)}=&
\mathbb{C}[SL(n)]^{N^+}&=^{\mathbb{C}[\delta^{(s)}_{i_1,\dots, i_s}]_{i_1>\dots
>i_s}}/_{P(\delta)}\\
\vspace{.1cm}\\
Vect(antisemistandard)=&&=Vect(semistandard)
\end{array}\\
\vspace{.1cm}\\
\swarrow~~~~~~~~~~\searrow\\
\vspace{.1cm}\\
\begin{array}{rcl}
^{\mathbb{C}[SL(n)]^{N^-}}/_{(P(\delta)+R(\delta)^-)}\simeq&\mathbb{C}[N^+]
\longleftrightarrow
\mathbb{C}[N^-]&\simeq^{\mathbb{C}[SL(n)]^{N^+}}/_{(P(\delta)+R(\delta)^+)}\\
\vspace{.1cm}\\
^{\mathbb{C}[\delta^{(s)}_{i_1, \dots,i_s}]_{i_s<n+1-s}}/_{P_{red}(\delta)^-}=&
&=^{\mathbb{C}[\delta^{(s)}_{i_1,\dots, i_s}]_{i_s>s}}/_{P_{red}(\delta)^-}\\
\vspace{.1cm}\\
Vect(antiquasistandard)=&&=Vect(quasistandard)
\end{array}\\
\end{array}$$

\

\section{The $\mathfrak{sl}(2)$ case}

\

\subsection{Representations of $\mathfrak{sl}(2)$}

\

The $\mathfrak{sl}(2)$-simple modules are charachterized by a
highest weight $a $. More precesely, the basis of
$\mathfrak{sl}(2)$ is:
$$X_\alpha=\left[\begin{matrix}0&1\\0&0\end{matrix}\right],\quad H_\alpha=\left[
\begin{matrix}1&0\\0&-1\end{matrix}\right],\quad Y_\alpha=\left[\begin{matrix}0&
0\\1&0\end{matrix}\right].$$ If $a$ is a positive integer, the
simple module $\pi^a$ acting on the space $V^a$ is
$a+1$-dimensional, with a basis $v_n$ ($0\leq n\leq a$) and the
matrices of the action are:
\begin{align*}
\pi^a(X_\alpha)&=\left[\begin{matrix}0&1&0&&0\\0&0&2&&0\\&&\ddots&\ddots&\\
&&&&a\\0&&&&0\end{matrix}\right],\\ \pi^a(H_\alpha)&=\left[
\begin{matrix}a&0&&0\\0&a-2&&0\\&&\ddots&\\0&&&-a\end{matrix}
\right],\\
\pi^a(Y_\alpha)&=\left[\begin{matrix}0&0&&0&0\\a&0&&0&0\\&&\ddots&\ddots&\\
0&&&1&0\end{matrix}\right].
\end{align*}

\

There is only one fundamental representation, associated to the
weight $\omega_1 $. We realize it in the space generated by the
functions $\delta^{(1)}_1(g) =g_{11}$, $\delta^{(1)}_2(g)=g_{21}$.
The other representations are realized on the space of homogeneous
polynomial functions of degree $a$ in these variables.

\

\subsection{Shape and reduced shape algebra}

\

There are no Pl\"ucker relation between $g_{11}$ and $g_{22}$,
thus the shape algebra is isomorphic to the algebra
$$A^\bullet(V)={\mathbb C}[g_{11},g_{21}]\simeq{\mathbb S}(V).$$
The reduced shape algebra is the quotient by the ideal generated
by $g_{11}-1$. Let us put:
$$\mathfrak{n}^-=\left\{\left[\begin{matrix}0&0\\x&0\end{matrix}\right]\right\}
,\quad
N^-=exp\left(\mathfrak{n}^-\right)=\left\{\left[\begin{matrix}1&0\\x&1
\end{matrix}\right]\right\}.$$
Then:
$${\mathbb S}_{red}^\bullet(V)^+={\mathbb C}[\delta^1_2]={\mathbb C}[X],$$
The $X_\alpha$ acts on a polynomial function as the operator:
$$X_\alpha=\frac{\partial}{\partial X}.$$
We realize the $\mathfrak{sl}(2)$-diamond cone as the half line of
the entire nodes $0,1,\dots,a,a+1,\dots$, at each node $n$, we put
the quasi standard Young tableau $\begin{array}{|c|c|c|} \hline
2&\dots & 2 \\
\hline\end{array}$ or the monomial $X^n$. We have an explicit
basis for the representation of $N^+$ on the diamond cone defined
by the action of $X_\alpha$, pictured by the graph:

\

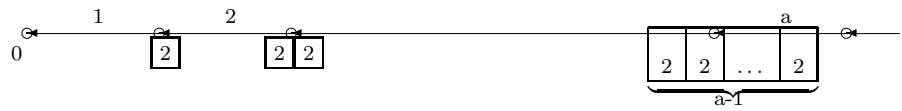
\begin{center}
{\tiny
\begin{picture}(410,600)(0,0)
\path(370,250)(40,250)

\put(92,250){\vector(-1,0){2}} \put(142,250){\vector(-1,0){2}}
\put(42,250){\vector(-1,0){2}} \put(352,250){\vector(-1,0){2}}
\put(302,250){\vector(-1,0){2}}

\put(40,250){\circle{4}}
\put(90,250){\circle{4}} \put(140,250){\circle{4}}
\put(300,250){\circle{4}} \put(350,250){\circle{4}}

\put(84,240){ \framebox{2}}
\put(130,240){\framebox{2}\framebox{2}} \put(34,240){0}
\put(275,240){$\underbrace{\begin{array}{|c|c|c|c|} \hline &&&\\
 2 & 2 & \dots
& 2
\\\hline
\end{array}}$}
\put(65,254){1} \put(115,254){2}\put(325,254){a}
\put(300,223){a-1}

\end{picture}}

\end{center}

\

For any $a\geq0$, we define the diamond $D_a$ as the graph
generated by $X^a$, the vector space $V^a$ as the vector space
with basis the nodes of $D_a$.\\
We saw that the anti semi standard (resp. the anti quasi standard)
basis can be identified with the semi standard (resp. the quasi
standard) basis. More precisely, a being fixed, the action of
$\tau$ on $V^a$, denoted by $\tau^{(a)}$ is defined as:

\

$$ \tau^{(a)}(X^n)=X^{a-n}$$

$$\tau^{(a)}( \begin{array}{|c|c|c|} \hline 2 & \dots &2 \\ \hline \end{array})
=\begin{array}{|c|c|c|} \hline 2 & \dots &2 \\ \hline \end{array}
$$

\

We can see $\tau^{(a)}$ as the succession of the operations:
\begin{itemize}
\item Completion of the tableau $T$ ($compl[2\dots 2]=[1\dots
1][2\dots 2]$) \item Action of $\tau$ ($\tau (compl[2\dots
2])=[2\dots 2][1\dots 1]$) \item reordering ($ord (\tau(
compl[2\dots 2]))=[1\dots 1][2\dots 2]$) \item Cancelling the
trivial columns $\begin{array}{|c|} \hline 1\\ \hline \end{array}$
\end{itemize}
We put:
$$Y_\alpha(X^n)=(\tau^{(a)}\circ X_\alpha\circ\tau^{(a)})(X^n)=(a-n)X^{n+1}$$
and $H_\alpha=[X_\alpha,Y_\alpha]$ or:
$$H_\alpha(X^n)=[(n+1)(a-n)-n(a-n+1)]X^n=(a-2n)X^n.$$
We complete the diamond $D_a$ by adding the edges corresponding to
the $Y_\alpha $-action.

\

\

\section{The $\mathfrak{sl}(3)$ case}

\

\subsection{Representations of $\mathfrak{sl}(3)$}

\

The $\mathfrak{sl}(3)$-simple modules are characterized by their
highest weight. More precisely, the basis of $\mathfrak{sl}(3)$
is:
\begin{align*}
X_\alpha&=\left[\begin{matrix}0&1&0\\0&0&0\\0&0&0\end{matrix}\right],
\quad
X_\beta=\left[\begin{matrix}0&0&0\\0&0&1\\0&0&0\end{matrix}\right],
\quad
X_{\alpha+\beta}=\left[\begin{matrix}0&0&1\\0&0&0\\0&0&0\end{matrix}
\right],\\
H_\alpha&=\left[\begin{matrix}1&0&0\\0&-1&0\\0&0&0\end{matrix}\right],
\quad H_\beta=\left[\begin{matrix}0&0&0\\0&1&0\\0&0&-1\end{matrix}\right],\\
Y_\alpha&=\left[\begin{matrix}0&0&0\\1&0&0\\0&0&0\end{matrix}\right],
\quad
Y_\beta=\left[\begin{matrix}0&0&0\\0&0&0\\0&1&0\end{matrix}\right],
\quad
Y_{\alpha+\beta}=\left[\begin{matrix}0&0&0\\0&0&0\\1&0&0\end{matrix}
\right].
\end{align*}
The simple modules have non multiplicity free weights. We can
describe then by using the reduced shape algebra. The fundamental
modules are three dimensional, they are realized on the space
$V^{\omega_1}={\mathbb C}^3$ and $V^{\omega_2}= \wedge^2{\mathbb
C}^3$.

\noindent For each pair of natural integers, there is an unique
irreducible representation $\pi(a,b)$ with highest weight
$a\varpi_1 + b\varpi_2$.

 \

\subsection{Shape and reduced shape algebra}

\

Now we have just one Pl\"ucker relation: let us put as above:
\begin{align*}
&\delta^{(1)}_1=g_{11},\quad\quad&\delta^{(1)}_2=g_{21},\quad\quad\quad\quad
\quad\quad&\delta^{(1)}_3=g_{31}\hfill\\
&\delta^{(2)}_{12}=g_{11}g_{22}-g_{12}g_{21},\quad&\delta^{(2)}_{13}=g_{11}
g_{32}-g_{12}g_{31},\quad&\delta^{(2)}_{23}=g_{21}g_{32}-g_{22}g_{31}.
\end{align*}
Then the unique Pl\"ucker relation is:
$$\delta^{(1)}_1\delta^{(2)}_{23}-\delta^{(1)}_2\delta^{(2)}_{13}+\delta^{(1)}_3
\delta^{(2)}_{12}=0.$$ The shape algebra is the quotient of the
algebra of polynomial functions in these 6 variables by the above
relation.

The reduced shape algebra is obtained by imposing
$\delta^{(1)}_1=1$ and $\delta^{(2)}_{12}=1$.

An explicit description of a basis for this module $V^{(a,b)}$ and
the $X_{\eta}, ~Y_{\eta},~H~_{\eta}$ actions on this basis can be
found in \cite{W} for instance. More precisely, Wildberger defines a
diamond cone $D$ in $\mathbb{R}^3$ and a infinite dimensional
vector space $V$ with basis: $$  \begin{array}{cl} \mathcal{B} &=
\{ e_{m,n,\ell},~ (m,n,\ell)\in
D \subset \mathbb{R}^3 \}\\
&=\{e_{m,n,\ell},~ m,n \geq 0, ~ -n \leq \ell \leq 2m-n,~m-2n\leq\ell\leq m,\cr
&\hskip 7cm~ \ell\equiv max(m,n)mod 2\}.
\end{array}$$
He defines the action of $X_{\eta}$ on these vectors
$e_{m,n,\ell}$ and the irreducible module $V^{(a,b)}$ with highest
weight $a \varpi_1 + b \varpi_2$ is the module generated by the
$X_{\eta}$ action on the highest weight vector $e_{a+b, a+b,
a-b}$.

A basis for this module is an explicit subset
$\mathcal{B}^{(a,b)}$ of $\mathcal{B}$. There is a symmetry
$\tau_{(a,b)}$ on $V^{(a,b)}$, $\tau_{(a,b)}(\mathcal{B}^{(a,b)})
= \mathcal{B}^{(a,b)}$ and the $Y_{\eta},~H_{\eta}$ actions are
defined as: $$ Y_{\eta} = \tau_{(a,b)}\circ X_{\eta} \circ
\tau_{(a,b)},~~~H_{\eta} = [X_{\eta}, Y_{\eta}]$$ see \cite{W} for
explicit formulas.

\

Let us put:
$$\mathfrak{n}^-=\left\{\left[\begin{matrix}0&0&0\\x&0&0\\u&y&0\end{matrix}
\right]\right\}
,\quad N^-=exp\left(\mathfrak{n}^-\right)=\left\{\left[\begin{matrix}1&0&0\\
x&1&0\\u+\frac{xy}{2}&y&1\end{matrix}\right]\right\}.$$ Then:
$$\delta^{(1)}_2=X,\quad\delta^{(1)}_3=\frac{xy}{2}+u=U,\quad
\delta^{(2)}_{13}=Y,\quad\delta^{(2)}_{23}=\frac{xy}{2}-u=E$$ and
\begin{align*}
{\mathbb S}_{red}^\bullet(V)^+&\simeq{\mathbb C}[x,y,u]\\
&=~^{\hbox{${\mathbb
C}[\delta^{(1)}_2,\delta^{(1)}_3,\delta^{(2)}_{13},
\delta^{(2)}_{23}]$}}/_{\hbox{$\langle\delta^{(1)}_3+\delta^{(2)}_{23}-
\delta^{(1)}_2\delta^{(2)}_{13}\rangle$}}\\
&=~^{\hbox{${\mathbb C}[X,Y,U,E]$}}/_{\hbox{$\langle U+E- XY
\rangle$}}.
\end{align*}

The quasi standard ordering on variables is:
$$\delta^{(1)}_3<\delta^{(1)}_2<\delta^{(2)}_{23}<\delta^{(2)}_{13},
~~~~~~~~~~ {\rm or}~~~~~~~~ U>X>E>Y.$$ Then the leading term for
this basis is $\delta^{(1)}_2 \delta^{(2)}_{13}=XY$, thus we get
the basis:
\begin{align*}
&\left\{(\delta^{(1)}_3)^u(\delta^{(2)}_{23})^e(\delta^{(1)}_2)^x=U^u
E^eX^x,\quad u,~e,~x\in{\mathbb N}\right\}\\
&\hskip 2cm\bigcup\left\{(\delta^{(1)}_3)^u(\delta^{(2)}_{23})^e
(\delta^{(1)}_2)^y=U^uE^eY^y,\quad u,~e,~y\in{\mathbb N},~ y
>0 \right\}.
\end{align*}

\

Now the action of $X_\alpha$, $X_\beta$ and $X_{\alpha+\beta}$ on
these polynomials are the following:
$$X_\alpha=\frac{\partial}{\partial x}-\frac{y}{2}\frac{\partial}{\partial u},
\quad X_\beta=\frac{\partial}{\partial
y}+\frac{x}{2}\frac{\partial}{\partial u} ,\quad
X_{\alpha+\beta}=\frac{\partial}{\partial u},$$ or
\begin{align*}
&X_\alpha(X)=1,\quad&X_\alpha(Y)=0,\quad&X_\alpha(U)=0,\quad&
X_\alpha(E)=Y,\quad\\
&X_\beta(X)=0,\quad&X_\beta(Y)=1,\quad&X_\beta(U)=X,\quad&X_\beta
(E)=0,\quad\\
&X_{\alpha+\beta}(X)=0,\quad&X_{\alpha+\beta}(Y)=0,\quad&X_{\alpha+\beta}
(U)=1,\quad&X_{\alpha+\beta}(E)=-1.
\end{align*}
Then the $X_\eta$ acting by derivations on the polynomial
functions $\phi$, we refind the diamond cone, the diamond
$D^{(a,b)}$, the vector space $V^{(a,b)}$, the symmetry
$\tau_{(a,b)}$ and the complete diamond graphs on $D^{(a,b)}$
described in \cite{W} with the identification:
$$
\begin{array}{ll}
e_{m,n,\ell} = U^{n-\frac{m-\ell}{2}}
E^{\frac{m-\ell}{2}} X^{m-n} ~~~&{\rm if}~~~m > n\\
e_{m,n,\ell} = U^{\frac{m+\ell}{2}} E^{\frac{m-\ell}{2}} ~~~&{\rm if}~~~m = n\\
e_{m,n,\ell} = U^{\frac{m+\ell}{2}} E^{m-\frac{n+\ell}{2}} Y^{n-m}
~~~&{\rm if}~~~m < n,
\end{array}
$$
our basis coincide with the basis $\mathcal{B}$ given by Wildberger, for 
$\mathbb{S}_{red}^{\bullet}(V)^+$.\\

\

\subsection{$X_\eta$ action, Symmetry and $Y_\eta$ action}

\

With our notations, we have the following identification between
column and variables $X, U, Y, E$:
\begin{displaymath}
\begin{array}{ccc}
X=\delta^{(1)}_2(g)&~~\longrightarrow~~&\begin{array}{|c|}\hline 2\\
\hline\end{array}\cr &&\cr
U=\delta^{(1)}_3(g)&~~\longrightarrow~~&\begin{array}{|c|}\hline 3\\
\hline\end{array}\cr &&\cr
Y=\delta^{(2)}_{13}(g)&~~\longrightarrow~~&\begin{array}{|c|}\hline 1\\
\hline 3\\ \hline
\end{array}\cr
&&\cr
E=\delta^{(2)}_{23}(g)&~~\longrightarrow~~&\begin{array}{|c|}\hline 2\\
\hline 3\\ \hline\end{array}\cr
\end{array}\end{displaymath}
The unique reduced Pl\"ucker relation is:
\begin{displaymath}
\begin{array}{ccccccc}
\begin{array}{c}
\begin{array}{|c|}
\hline
3\\
\hline
\end{array}\\
\end{array}
&-&
\begin{array}{c}
\begin{array}{|c|c|}
\hline
1&2\\
\hline
\end{array}\\
\begin{array}{|c|}
3\\
\hline
\end{array}\hfill
\end{array}&+&
\begin{array}{c}
\begin{array}{|c|}
\hline
2\\
\hline
\end{array}\\
\begin{array}{|c|}
3\\
\hline\end{array}\hfill
\end{array}&=&0
\end{array}
\end{displaymath}

For instance the $X_{\alpha}$ action on our basis is exactly the Wildberger's 
one:

\

\begin{displaymath}\begin{array}{|c|c|}
  \hline
  U^u E^e X^x & e U^{u+1} E^{e - 1} X^{x-1} + (e + x)U^u E^e X^{x - 1}~~(x>0)\\
  \hline
  {U^u E^e Y^y} & e U^u E^{e - 1} Y^{y+1} ~~ (y \geq 0) \\ \hline
\end{array}
\end{displaymath}

\

or
\begin{displaymath}\begin{array}{|c|c|}
  \hline
  e_{m,n,\ell} & \frac{m-\ell}{2} e_{m-1,n,\ell+1} + (m-n+\frac{m-\ell}{2})
  e_{m-1, n, \ell-1} ~~ (m>n) \\
  \hline
  e_{m,n,\ell} & (m-\frac{n+\ell}{2}) e_{m-1,n,\ell} ~~ (n \geq m) \\ \hline
\end{array}
\end{displaymath}

\

\noindent and the $X_{\beta}$ action:

\

\begin{displaymath}
\begin{array}{|c|c|}
  \hline
  U^u E^e X^x & u U^{u-1} E^{e} X^{x-1}  ~~ (x>0) \\
  \hline
  U^u E^e Y^y & u U^{u-1} E^{e} X Y^{y} + y U^u E^e Y^{y-1} ~~ (y \geq 0) \\ \hline
\end{array}
\end{displaymath}

\

\noindent
or
\begin{displaymath}\begin{array}{|c|c|}
  \hline
  e_{m,n,\ell} & (n-\frac{m-\ell}{2}) e_{m,n-1,\ell}  ~~ (m>n) \\
  \hline
  e_{m,n,\ell} & (n-m+\frac{n+\ell}{2}) e_{m,n-1,\ell+1} + (\frac{n+\ell}{2})
e_{m,n-1,\ell-1} ~~ (n \geq m) \\ \hline
\end{array}
\end{displaymath}

\

For $\frak{sl}(3)$, our symmetry $\tau$ on quasi standard Young
tableaux induces a very simple transformation on $V^{(a,b)}$.

Starting with a quasi standard Young tableau $T$ with $a'$ columns
of height $1$ and $b'$ columns of height $2$, $a'\leq a$ and
$b'\leq b$, we complete $T$ by adding $a-a'$ trivial columns
$\begin{array}{|c|} \hline 1\\ \hline \end{array}$ and $b'-b$
trivial columns $\begin{array}{|c|} \hline 1\\
\hline 2\\
\hline \end{array}$. For instance:
$$a=5,~~ b=3~~~~~~~T=\begin{array}{l}
\begin{array}{|c|c|c|c|c|}
\hline 2&2&2&2&3\\
\hline  \end{array}\\
\begin{array}{|c|c|}
\hline 3 & 3\\
\hline
\end{array}\\
\end{array}, ~~ compl(T)=\begin{array}{l}
\begin{array}{|c|c|c|c|c|c|c|c|}
\hline 1&2&2&1&1&2&2&3\\
\hline  \end{array}\\
\begin{array}{|c|c|c|}
\hline 2 & 3 &3\\
\hline
\end{array}\\
\end{array}$$
Then we compute $\tau(compl(T))$, we reorganize the columns as
above and finally we suppress the trivial columns, on our example:
$$\begin{array}{ccl}
\tau (compl(T))&=& \begin{array}{l}
\begin{array}{|c|c|c|c|c|c|c|c|}
\hline 3&2&2&3&3&2&2&1\\
\hline  \end{array}\\
\begin{array}{|c|c|c|}
\hline 2 & 1 &1\\
\hline
\end{array}\\
\end{array}\\
&&\\
&=& - \begin{array}{l}
\begin{array}{|c|c|c|c|c|c|c|c|}
\hline 1&1&2&1&2&2&3&3\\
\hline  \end{array}\\
\begin{array}{|c|c|c|}
\hline 2 & 2 &3\\
\hline
\end{array}\\
\end{array}\\
&&\\
&\simeq& -\begin{array}{l}
\begin{array}{|c|c|c|c|c|}
\hline 2&2&2&3&3\\
\hline  \end{array}\\
\begin{array}{|c|}
\hline 3\\
\hline
\end{array}\\
\end{array}
\end{array}$$
The resulting quasi standard tableau will be denoted
$\tau^{(a,b)}(T)$. Explicitly we get with the polynomial
notations:
$$\begin{array}{ccc} \tau^{(a,b)} (U^u E^e X^x)& =&
U^{a-(x+u)} E^{b-e} X^x\\
\tau^{(a,b)} (U^u E^e X^x)& =& U^{a-u} E^{b-(y+e)} Y^y \end{array}
$$
or: $$ \tau^{(a,b)}(e_{m,n,\ell}) = e_{a+b-n, a+b-m,
a-b+m-n-\ell}$$ We refind the symmetry, thus the $Y_{\eta}$ and
$H_{\eta}$ actions of Wildberger (\cite{W}).

\

\section{The $\mathfrak{sl}(4)$ case}

\

\subsection{Representations of $\mathfrak{sl}(4)$}

\

As above, we have simple roots $\alpha$, $\beta$ and $\gamma$,
with:
$$X_\alpha=\left[\begin{matrix}0&1&0&0\\0&0&0&0\\0&0&0&0\\0&0&0&0\end{matrix}
\right],\quad
X_\beta=\left[\begin{matrix}0&0&0&0\\0&0&1&0\\0&0&0&0\\0&0&0&0
\end{matrix}\right],\quad X_\gamma=\left[\begin{matrix}0&0&0&0\\0&0&0&0\\
0&0&0&1\\0&0&0&0\end{matrix}\right].$$ Moreover we have positive
roots $\alpha+\beta$, $\beta+\gamma$ and $\alpha+\beta +\gamma$,
with:
$$X_{\alpha+\beta}=\left[\begin{matrix}0&0&1&0\\0&0&0&0\\0&0&0&0\\0&0&0&0
\end{matrix}\right],\quad X_{\beta+\gamma}=\left[\begin{matrix}0&0&0&0\\
0&0&0&1\\0&0&0&0\\0&0&0&0\end{matrix}\right],\quad
X_{\alpha+\beta+\gamma}=
\left[\begin{matrix}0&0&0&1\\0&0&0&0\\0&0&0&0\\0&0&0&0\end{matrix}\right].$$
We put $Y_\eta=~^tX_\eta$ and
$$H_\alpha=\left[\begin{matrix}1&0&0&0\\0&-1&0&0\\0&0&0&0\\0&0&0&0\end{matrix}
\right],\quad
H_\beta=\left[\begin{matrix}0&0&0&0\\0&1&0&0\\0&0&-1&0\\0&0&0&0
\end{matrix}\right],\quad H_\gamma=\left[\begin{matrix}0&0&0&0\\0&0&0&0\\
0&0&1&0\\0&0&0&-1\end{matrix}\right].$$ The fundamental
representations are 4 and 6 dimensional, they are associted to the
fundamental highest weight $\omega_1$ for the canonical
representation on $V= {\mathbb C}^4$, $\omega_2$ for the
representation on $\wedge^2V$ and $\omega_3$ for the
representation on $\wedge^3V$. These fundamental representations
are easy to describe, the reduction of the tensor product of any
two of them is completely described in \cite{FH}.
Especially, we get the Pl\"ucker relations via this decomposition.

\

\subsection{Shape and reduced shape algebra}

\

Now we have 10 Pl\"ucker relations: let us put as above:
$$\delta^{(1)}_i=g_{i1},\quad\delta^{(2)}_{ij}=\left|\begin{matrix}g_{i1}&g_{i2}
\\g_{j1}&g_{j2}\end{matrix}\right|,\quad
\delta^{(3)}_{ijk}=\left|\begin{matrix}g_{i1}&g_{i2}&g_{i3}\\g_{j1}&g_{j2}&
g_{j3}\\g_{k1}&g_{k2}&g_{k3}\end{matrix}\right|$$ Then we have 4
Pl\"ucker relations between the $\delta^{(1)}$ and $\delta^{(2)}$
:
\begin{align*}
\delta^{(1)}_1\delta^{(2)}_{23}-\delta^{(1)}_2\delta^{(2)}_{13}+\delta^{(1)}_3
\delta^{(2)}_{12}&=0,\\
\delta^{(1)}_2\delta^{(2)}_{34}-\delta^{(1)}_3\delta^{(2)}_{24}+\delta^{(1)}_4
\delta^{(2)}_{23}&=0,\\
\delta^{(1)}_1\delta^{(2)}_{34}-\delta^{(1)}_3\delta^{(2)}_{14}+\delta^{(1)}_4
\delta^{(2)}_{13}&=0,\\
\delta^{(1)}_1\delta^{(2)}_{24}-\delta^{(1)}_2\delta^{(2)}_{14}+\delta^{(1)}_4
\delta^{(2)}_{12}&=0.
\end{align*}
There are also 4 relations between the $\delta^{(2)}$ and
$\delta^{(3)}$:
\begin{align*}
\delta^{(2)}_{14}\delta^{(3)}_{234}-\delta^{(2)}_{24}\delta^{(3)}_{134}+
\delta^{(2)}_{34}\delta^{(3)}_{124}&=0,\\
\delta^{(2)}_{12}\delta^{(3)}_{134}-\delta^{(2)}_{13}\delta^{(3)}_{124}+
\delta^{(2)}_{14}\delta^{(3)}_{123}&=0,\\
\delta^{(2)}_{12}\delta^{(3)}_{234}-\delta^{(2)}_{23}\delta^{(3)}_{124}+
\delta^{(2)}_{24}\delta^{(3)}_{123}&=0,\\
\delta^{(2)}_{13}\delta^{(3)}_{234}-\delta^{(2)}_{23}\delta^{(3)}_{134}+
\delta^{(2)}_{34}\delta^{(3)}_{123}&=0.
\end{align*}
And one between the $\delta^{(2)}$:
$$\delta^{(2)}_{12}\delta^{(2)}_{34}-\delta^{(2)}_{13}\delta^{(2)}_{24}+
\delta^{(2)}_{14}\delta^{(2)}_{23}=0.$$ And finally one between
the $\delta^{(1)}$ and the $\delta^{(3)}$:
$$\delta^{(1)}_1\delta^{(3)}_{234}-\delta^{(1)}_2\delta^{(3)}_{134}+
\delta^{(1)}_3\delta^{(3)}_{124}-\delta^{(1)}_4\delta^{(3)}_{123}=0.$$

The shape algebra is the quotient of the algebra of polynomial
functions in these 14 variables by the 10 above relations.

The reduced shape algebra is obtained by imposing
$\delta^{(1)}_1=1$, $\delta^{(2)}_{12}=1$ and
$\delta^{(3)}_{123}=1$.

\

Let us put:
$$\mathfrak{n}^-=\left\{\left[\begin{matrix}0&0&0&0\\x&0&0&0\\u&y&0&0\\w&v&z&0
\end{matrix}\right]\right\}$$
and
$$N^-=exp(\mathfrak{n}^-)=\left\{\left[\begin{matrix}
1&0&0&0\\x&1&0&0\\u+\frac{xy}{2}&y&1&0\\w+\frac{xv}{2}+\frac{zu}{2}+
\frac{xyz}{6}&v+\frac{yz}{2}&z&1\end{matrix}\right]\right\}.$$
Then we get:
$$\delta^{(1)}_1=1,~~\delta^{(1)}_2=X,~~\delta^{(1)}_3=\frac{xy}{2}+u=U,
~~\delta^{(1)}_4=w+\frac{xv}{2}+\frac{zu}{2}+\frac{xyz}{6}=A$$ and
$$\delta^{(3)}_{123}=1,~~\delta^{(3)}_{124}=Z,~~\delta^{(3)}_{134}=yz-\xi_2=
W,~~\delta^{(3)}_{234}=\frac{xyz}{6}-\frac{xv}{2}-\frac{zu}{2}+w=C$$
and
\begin{displaymath}
\begin{array}{ccc}
\delta^{(2)}_{12}=1,\hfill&\delta^{(2)}_{13}=Y,\hfill&\delta^{(2)}_{14}=v+
\displaystyle{\frac{yz}{2}}=V\hfill\\
\delta^{(2)}_{23}=\displaystyle{\frac{xy}{2}-u=E},~~&\delta^{(2)}_{24}=
\displaystyle{\frac{xyz}{3}+\frac{xv}{2}-\frac{zu}{2}-w=D},~~&\delta^{(2)}_{34}=
\displaystyle{\frac{xy^2z}{12}+uv-yw=B.}
\end{array}
\end{displaymath}
Now:
\begin{align*}
{\mathbb S}^\bullet_{red}(V)^+&\simeq{\mathbb C}[x,y,z,u,v,w]\\
&=~~^{\hbox{${\mathbb
C}[\delta^{(1)}_2,\dots,\delta^{(1)}_4,\delta^{(2)}_{13},
\dots,\delta^{(2)}_{34},\delta^{(3)}_{124},\dots,\delta^{(3)}_{234}]$}}/_{\hbox
{$P_{red}(\delta)^+$}}\\
&=~~^{\hbox{${\mathbb C}[X,Y,Z,U,E,W,V,A,C,D,B]$}}/_{\hbox
{$Pluck$}}
\end{align*}
where $Pluck$ is the ideal generated by the 10 polynomials:
\begin{align*}
Pluck=\langle&U-XY+E,~~D-XV+A,~~B-UV+YA,~~XB-UD
+AE,\\
&B-YD+EV,~~C-XW_1+UZ-A,\\
&VC-DW+PZ,~~W-YZ+V,~~C-EZ+D,~~YC-EW+B \rangle.
\end{align*}

\

We choosed the following ordering for our variables:
$$
Z<W<C<Y<E<V<D<B<X<U<A.
$$
Then the leading terms of this basis are:
$$XY,~~XV,~~UV,~~BX,~~YZ,~~EZ,~~YC,~~VC,~~XW,~~EV,~~UDW,~~UDY.$$
Now the basis of our space, {\sl i.e.} the nodes of the
$\mathfrak{sl}(4)$-diamond are monomials
$$
X^x Y^y Z^z W^w V^v U^u E^e A^a C^c D^d B^b
$$
with:
$$
0=xy=xv=uv=bx=yz=ez=yc=vc=xw=ev=udw=udy.
$$

\

The action of our generators $X_\alpha$, $X_\beta$ and $X_\gamma$
on these polynomials are:
\begin{align*}
X_\alpha&=\partial_x-\frac{y}{2}\partial_u+\left(\frac{yz}{12}-\frac{v}{2}
\right)\partial_w,\\
X_\beta&=\partial_y+\frac{x}{2}\partial_u-\frac{z}{2}\partial_v-\frac{xz}{6}
\partial_w,\\
X_\gamma&=\partial_z+\frac{y}{2}\partial_v+\left(\frac{xy}{12}+\frac{u}{2}
\right)\partial_w.
\end{align*}
Then we get:
\begin{displaymath}
\displaystyle
\begin{array}{ccc}
&&\\
X_\alpha(X)=1,\hfill&X_\beta(X)=0,\hfill&X_\gamma(X)=0,\\&&\\
X_\alpha(Y)=0,\hfill&X_\beta(Y)=1,\hfill&X_\gamma(Y)=0,\\&&\\
X_\alpha(Z)=0,\hfill&X_\beta(Z)=0,\hfill&X_\gamma(Z)=1,\\&&\\
X_\alpha(U)=0,\hfill&X_\beta(U)=X,\hfill&X_\gamma(U)=0,\\&&\\
X_\alpha(E)=Y,\hfill&X_\beta(E)=0,\hfill&X_\gamma(E)=0,\\&&\\
X_\alpha(W)=0,\hfill&X_\beta(W)=Z,\hfill&X_\gamma(W)=0,\\&&\\
X_\alpha(V)=0,\hfill&X_\beta(V)=0,\hfill&X_\gamma(V)=Y,\\&&\\
X_\alpha(A)=0,\hfill&X_\beta(A)=0,\hfill&X_\gamma(A)=U,\\&&\\
X_\alpha(C)=W,\hfill&X_\beta(C)=0,\hfill&X_\gamma(C)=0,\\&&\\
X_\alpha(D)=V,\hfill&X_\beta(D)=0,\hfill&X_\gamma(D)=V,\\&&\\
X_\alpha(B)=0,\hfill&X_\beta(B)=D,\hfill&X_\gamma(B)=0.\\&&
\end{array}
\end{displaymath}

\

Thus the $X_\eta$ for $\eta$ simple are acting on our basis of the
reduced shape algebra by giving linear combination with integral
coefficients, indeed, we find first such a linear combination on
${\mathbb Z}$ (even ${\mathbb Z}^+$) coefficients but on monomials
which are perhaps not all admissible, then we come back to
admissible monomials, using the reduced Pl\"ucker relations, but
these relations are with coefficients $\pm1$, thus we finally get
a combination of monomials in the basis with coefficients in
${\mathbb Z}$.

\

\subsection{Symmetry}

\

Now the symmetry $\tau$ on Young tableaux does not induce a simple
operation $\tau^{(abc)}$ on the basis of the simple module
$V^{(abc)}$.\\
For instance the tableau $\begin{array}{l}
\begin{array}{|c|c|}
\hline 1&3\\
\hline
\end{array}\\
\begin{array}{|c|}
\hline 2\\
\hline
\end{array}\\
\begin{array}{|c|}
\hline 4\\
\hline
\end{array}\\
\end{array}=ZU$ is an element of the basis of $V^{(1,0,1)}$ (see
fig.). Repeating the operation performed for $\frak{sl}(2)$ and
$\frak{sl}(3)$, we get:
$$compl(\begin{array}{l}
\begin{array}{|c|c|}
\hline 1&3\\
\hline
\end{array}\\
\begin{array}{|c|}
\hline 2\\
\hline
\end{array}\\
\begin{array}{|c|}
\hline 4\\
\hline
\end{array}\\
\end{array})=\begin{array}{l}
\begin{array}{|c|c|}
\hline 1&3\\
\hline
\end{array}\\
\begin{array}{|c|}
\hline 2\\
\hline
\end{array}\\
\begin{array}{|c|}
\hline 4\\
\hline
\end{array}\\
\end{array},~~~~~ \tau(compl(\begin{array}{l}
\begin{array}{|c|c|}
\hline 1&3\\
\hline
\end{array}\\
\begin{array}{|c|}
\hline 2\\
\hline
\end{array}\\
\begin{array}{|c|}
\hline 4\\
\hline
\end{array}\\
\end{array}))=\begin{array}{l}
\begin{array}{|c|c|}
\hline 4&2\\
\hline
\end{array}\\
\begin{array}{|c|}
\hline 3\\
\hline
\end{array}\\
\begin{array}{|c|}
\hline 1\\
\hline
\end{array}\\
\end{array}=-\begin{array}{l}
\begin{array}{|c|c|}
\hline 1&2\\
\hline
\end{array}\\
\begin{array}{|c|}
\hline 3\\
\hline
\end{array}\\
\begin{array}{|c|}
\hline 4\\
\hline
\end{array}\\
\end{array}$$
but this tableau is not quasi standard: the extraction of the
trivial top $\begin{array}{|c|} \hline 1\\
\hline \end{array}$ of the first column is not trivial:
$$\tau^{(1,0,1)}(\begin{array}{l}
\begin{array}{|c|c|}
\hline 1&3\\
\hline
\end{array}\\
\begin{array}{|c|}
\hline 2\\
\hline
\end{array}\\
\begin{array}{|c|}
\hline 4\\
\hline
\end{array}\\
\end{array})=-\begin{array}{l}
\begin{array}{|c|}
\hline 2\\
\hline
\end{array}\\
\begin{array}{|c|}
\hline 3\\
\hline
\end{array}\\
\begin{array}{|c|}
\hline 4\\
\hline
\end{array}\\
\end{array}-\begin{array}{l}
\begin{array}{|c|c|}
\hline 1&3\\
\hline
\end{array}\\
\begin{array}{|c|}
\hline 2\\
\hline
\end{array}\\
\begin{array}{|c|}
\hline 4\\
\hline
\end{array}\\
\end{array}+\begin{array}{|c|}
\hline 4\\
\hline \end{array}$$ or:
$$\tau^{(1,0,1)}(ZU)=-WX=-C-ZU+A$$
We prefer to keep the new Groebner basis to see $\tau$ as a global
change of basis inside the reduced shape algebra and to realize
$Y_{\tau_\eta}=\tau X_\eta \tau$ by using the two basis. For
instance in $V^{(101)}$ the basis is:
$$\{1, X, U, A, Z, W, C, WU, WA, CU, CA, CX, ZU, ZA, ZX \}$$
the image by $\tau_{compl}$ of this basis is:
$$\{AC, UC, XC, C, WA, ZA, A, ZX, Z, X, 1, WX, W, WU\}$$
The matrix of $Y_{\tau_\eta}$ on this new basis is exactly the
matrix of $X_\eta$ in the old one.

Here is the presentation for adjoint representation $V^{(1,0,1)}$ of $SL(4)$:

\

{\tiny

\begin{picture}(410,600)(40,0)
\put(90,180){\vector(1,0){350}} \put(90,180){\vector(0,1){360}}
\put(90,180){\vector(-1,-2){90}} \path(90,180)(0,0)
\put(400,170){$(3,4)$} \put(95,530){$(2,3)$} \put(10,5){$(1,2)$}
\path(50,100)(50,280) \put(50,190){\vector(0,-1){1}}
\put(45,190){1} \path(50,100)(230,100)
\put(140,100){\vector(-1,0){1}} \put(140,90){1}
\dottedline[.](50,280)(220,280) \put(140,280){\vector(-1,0){1}}
\put(140,285){1} \dottedline[.](50,280)(230,270)
\put(140,275){\vector(-1,0){1}} \put(140,267){1}
\path(50,280)(230,460) \put(140,370){\vector(-1,-1){1}}
\put(135,375){-1} \path(50,280)(90,180)
\put(70,230){\vector(1,-2){1}} \put(65,225){1}
\path(50,280)(190,380) \put(120,330){\vector(-1,-1){1}}
\put(115,335){1} \path(190,380)(230,460)
\put(210,420){\vector(1,2){1}} \put(205,420){1}
\path(230,460)(410,460) \put(320,460){\vector(-1,0){1}}
\put(320,465){1} \path(230,460)(270,360)
\put(250,410){\vector(1,-2){1}} \put(250,415){1}
\dottedline[.](230,460)(240,290) \put(234,390){\vector(0,-1){1}}
\put(239,390){1} \dottedline[.](230,460)(220,280)
\put(226,390){\vector(0,-1){1}} \put(218,390){-1}
\dottedline[.](230,460)(230,270) \put(230,340){\vector(0,-1){1}}
\put(233,340){2} \path(410,460)(410,280)
\put(410,370){\vector(0,-1){1}} \put(415,370){1}
\path(410,460)(370,380) \put(390,420){\vector(1,2){1}}
\put(392,417){1} \path(410,460)(270,360)
\put(340,410){\vector(-1,-1){1}} \put(337,413){1}

\dottedline[.](410,460)(230,270)
\put(344,390){\vector(-1,-1){1}}
\put(350,390){1}

\dottedline[.](410,460)(220,280)
\put(336,390){\vector(-1,-1){1}}
\put(326,390){-2}

\dottedline[.](410,460)(240,290) \put(290,340){\vector(-1,-1){1}}

\put(292,345){1} \path(370,380)(190,380)
\put(290,380){\vector(-1,0){1}} \put(290,383){1}
\path(370,380)(410,280) \put(390,330){\vector(1,-2){1}}
\put(380,330){-1} \path(370,380)(190,200)
\put(300,310){\vector(-1,-1){1}} \put(306,310){1}

\dottedline[.](370,380)(220,280) \put(301,334){\vector(-1,-1){1}}
\put(298,336){1} \dottedline[.](370,380)(240,290)
\put(307,336){\vector(-1,-1){1}} \put(312,334){1}

\path(270,360)(90,180) \put(169,259){\vector(-1,-1){1}}
\put(162,262){-1} \path(270,360)(270,180)
\put(270,250){\vector(0,-1){1}} \put(274,251){1}
\dottedline[.](270,360)(240,290) \put(255,325){\vector(1,2){1}}
\put(260,325){1} \path(190,200)(50,100)
\put(120,150){\vector(-1,-1){1}} \put(126,148){1}
\path(190,200)(230,100) \put(210,150){\vector(1,-2){1}}
\put(214,150){-1} \path(190,200)(190,380)
\put(190,300){\vector(0,-1){1}} \put(184,300){1}
\dottedline[.](190,200)(240,290) \put(216,247){\vector(1,2){1}}
\put(210,245){2} \dottedline[.](190,200)(230,270)
\put(217,246){\vector(1,2){1}} \put(220,244){1}
\dottedline[.](190,200)(220,280) \put(205,240){\vector(1,2){1}}
\put(197,240){-1} \path(410,280)(270,180)
\put(340,230){\vector(-1,-1){1}} \put(335,232){1}
\path(410,280)(230,100) \put(320,190){\vector(-1,-1){1}}
\put(326,190){1}

\dottedline[.](410,280)(230,270) \put(320,275){\vector(-1,0){1}}
\put(320,267){1} \dottedline[.](410,280)(220,280)
\put(315,280){\vector(-1,0){1}} \put(315,284){1}
\path(270,180)(230,100) \put(250,140){\vector(1,2){1}}
\put(244,140){1} \dottedline[.](270,180)(240,290)
\put(255,235){\vector(1,-2){1}} \put(259,235){-1}
\dottedline[.](270,180)(230,270) \put(250,225){\vector(1,-2){1}}
\put(243,225){1} \dottedline[.](230,270)(230,100)
\put(230,185){\vector(0,-1){1}} \put(234,185){1}
\dottedline[.](220,280)(50,100) \put(135,190){\vector(-1,-1){1}}
\put(142,190){1} \dottedline[.](220,280)(90,180)
\put(155,230){\vector(-1,-1){1}} \put(153,234){1}
\dottedline[.](240,290)(90,180) \put(162,232){\vector(-1,-1){1}}
\put(165,228){1}

\dottedline[.](190,380)(240,290) \put(215,335){\vector(1,-2){1}}
\put(218,335){1}

\dottedline[.](190,380)(230,270) \put(208,330){\vector(1,-2){1}}
\put(200,330){-1}

\put(90,180){\circle{4}}
\put(50,100){\circle{4}}
\put(50,280){\circle{4}}
\put(230,460){\circle{4}}
\put(410,460){\circle{4}}
\put(410,280){\circle{4}}
\put(270,360){\circle{4}}
\put(190,380){\circle{4}}
\put(270,360){\circle{4}}
\put(190,200){\circle{4}}
\put(270,180){\circle{4}}
\put(230,100){\circle{4}}
\put(370,380){\circle{4}}
\put(220,280){\circle{4}}
\put(230,270){\circle{4}}
\put(240,290){\circle{4}}
\put(33,280){\framebox{3}}
\put(52,90){\framebox{2}}
\put(230,80){$\begin{array}{l}
\framebox{1}\framebox{2}\\
\framebox{2}\\
\framebox{4}\\
\end{array}$}
\put(271,161){$\begin{array}{l}
\framebox{1}\\
\framebox{2} \\
\framebox{4}\\
\end{array}$}
\put(193,190){$\begin{array}{l}
\framebox{2}\framebox{2}\\
\framebox{3} \\
\framebox{4}\\
\end{array}$}
\put(228,476){$\begin{array}{l}
\framebox{1}\framebox{3}\\
\framebox{3} \\
\framebox{4}\\
\end{array}$}
\put(409,476){$\begin{array}{l}
\framebox{1}\framebox{4}\\
\framebox{3} \\
\framebox{4}\\
\end{array}$}
\put(410,270){$\begin{array}{l}
\framebox{1}\framebox{4}\\
\framebox{2} \\
\framebox{4}\\
\end{array}$}
\put(168,390){$\begin{array}{l}
\framebox{2}\framebox{3}\\
\framebox{3} \\
\framebox{4}\\
\end{array}$}
\put(373,390){$\begin{array}{l}
\framebox{2}\framebox{4}\\
\framebox{3}\\
\framebox{4}\\
\end{array}$}
\put(271,348){$\begin{array}{l}
\framebox{1}\\
\framebox{3} \\
\framebox{4}\\
\end{array}$}
\put(246,282){$\begin{array}{l}
\framebox{2}\\
\framebox{3}\\
\framebox{4}\\
\end{array}$}
\put(218,248){$\begin{array}{l}
\framebox{1}\framebox{3}\\
\framebox{2}\\
\framebox{4}\\
\end{array}$}
\put(204,287){\framebox{4}}
\put(93,172){0}

\end{picture}}

\

\end{document}